\documentclass[10pt,a4paper,english,DIV=9]{scrartcl}
\pdfoutput=1

\usepackage[T1]{fontenc}
\usepackage[utf8]{inputenc}

\usepackage[sc]{mathpazo}   
\usepackage[scaled=.95]{helvet}
\usepackage{courier}

\usepackage{microtype}

\KOMAoptions{DIV=last}						

\usepackage{amssymb,amsthm, amsmath,mathtools}

\usepackage{dsfont}						

\usepackage{tikz}							
\usetikzlibrary[matrix,arrows]

\usepackage{subcaption}

\usepackage{tabularx}

\usepackage{hyperref}
\hypersetup{
  colorlinks   = true,          
  urlcolor     = blue,    
  linkcolor    = blue,       
  citecolor   = purple        
}

\let\OLDthebibliography\thebibliography
\renewcommand\thebibliography[1]{
  \OLDthebibliography{#1}
  \small
  \setlength{\parskip}{0pt}
  \setlength{\itemsep}{1.5pt plus 0.3ex}
}

 \theoremstyle{plain}
 \newtheorem{theorem}{Theorem}[section]
  \theoremstyle{definition}
  \newtheorem{definition}[theorem]{Definition}
  \theoremstyle{definition}
  \newtheorem{example}[theorem]{Example}
  \theoremstyle{plain}
  \newtheorem{lemma}[theorem]{Lemma}
  \theoremstyle{plain}
  \newtheorem{corollary}[theorem]{Corollary}
  \theoremstyle{remark}
  
  \theoremstyle{plain}
  \newtheorem{proposition}[theorem]{Proposition}
  \theoremstyle{plain}
  
  \theoremstyle{remark}
  
  \theoremstyle{remark}
  
  \theoremstyle{plain}
  
  \theoremstyle{definition}
  \newtheorem{construction}[theorem]{Construction}

\DeclareMathOperator{\Hom}{Hom}

\DeclareMathOperator{\Spec}{Spec}
\DeclareMathOperator{\Trop}{Trop}
\DeclareMathOperator{\trop}{trop}

\DeclareMathOperator{\Newt}{Newt}

\DeclareMathOperator{\ini}{in}

\DeclareMathOperator{\Pic}{Pic}
\DeclareMathOperator{\conv}{conv}
\DeclareMathOperator{\Lin}{Lin}

\DeclareMathOperator{\mult}{mult}

\DeclareMathOperator{\GL}{GL}
\DeclareMathOperator{\G}{G}
\DeclareMathOperator{\dir}{v}

\DeclareMathOperator{\dist}{dist}
\DeclareMathOperator{\Image}{Im}
\DeclareMathOperator{\pr}{pr}

\newcommand{\Mn}{M_{0,n}(r,d)}
\newcommand{\Mln}{M^{\operatorname{lab}}_{0,n}(r,d)}
\newcommand{\Mno}[1]{M_{0,#1}(r,d)}
\newcommand{\Mlno}[1]{M^{\operatorname{lab}}_{0,#1}(r,d)}
\newcommand{\extMlno}{\mathds G_m^{\dPair {L_0}}/\mathds G_m^{L_0}\times\aT r}
\newcommand{\MlIrd}[3]{M^{\operatorname{lab}}_{0,#1}(#2,#3)}

\newcommand{\dPair}[1]{{D(#1)}}
\newcommand{\VG}{\Gamma_{\nu}}

\newcommand{\PS}[1]{\ensuremath{\mathds P^{#1}}}
\newcommand{\aT}[1]{\ensuremath{\mathds T^{#1}}}

\newcommand{\tG}{\mathfrak G}
\newcommand{\tMno}[1]{\mathfrak M_{0,#1}}
\newcommand{\tMln}{\mathfrak M^{\operatorname{lab}}_{0,n}(r,d)}
\newcommand{\tMlI}[1]{\mathfrak M^{\operatorname{lab}}_{0,#1}(r,d)}

\newcommand{\tMlno}[3]{\mathfrak M^{\operatorname{lab}}_{0,#1}(#2,#3)}

\newcommand{\tpl}{\mathfrak{pl}}

\newcommand{\amin}{d}
\newcommand{\tmin}{\mathfrak d}

\newcommand{\apara}{f}
\newcommand{\tpara}{\mathfrak f}

\newcommand{\tT}[1]{\ensuremath{\boldsymbol{\mathfrak T}^{#1}}}
\newcommand{\ttorus}{\mathfrak T}
\newcommand{\tOrb}{\mathfrak O}

\newcommand{\TR}[1][{}]{\ensuremath{\mathcal T_{#1}(R)}}
\newcommand{\TpR}[1][{}]{\ensuremath{\mathcal T'_{#1}(R)}}
	
\newcommand{\ev}[1]{\operatorname{ev}_{#1}}
\newcommand{\tev}[1]{\ensuremath{\mathfrak{ev}_{#1}}}
\newcommand{\extev}[1]{\overline{\operatorname{ev}}_{#1}}
\newcommand{\textev}[1]{\ensuremath{\overline{\mathfrak{ev}}_{#1}}}

\newcommand{\wT}{\widetilde T}
\newcommand{\bT}{\overline T}
\newcommand{\wM}{\widetilde M}

\newcommand{\nr}{N_{\mathds R}}
\newcommand{\nt}{N_{\ttorus}}

\renewcommand{\phi}{\varphi}

\addtokomafont{section}{\large\rmfamily}
\addtokomafont{subsection}{\normalsize\rmfamily}
\addtokomafont{title}{\rmfamily\mdseries\scshape}	
\addtokomafont{paragraph}{\rmfamily}

\begin{document}

\title{\large Correspondence Theorems via Tropicalizations of Moduli Spaces}

\author{\normalsize Andreas Gross}
\date{}

\maketitle
\vspace{-1.2cm}

\begin{abstract}
We show that the moduli spaces of irreducible labeled parametrized marked rational curves in toric varieties can be embedded into algebraic tori such that their tropicalizations are the analogous tropical moduli spaces. These embeddings are shown to respect the evaluation morphisms in the sense that evaluation commutes with tropicalization. With this particular setting in mind we prove a general correspondence theorem for enumerative problems which are defined via ``evaluation maps'' in both the algebraic and tropical world. Applying this to our motivational example we reprove Nishinou and Siebert's correspondence theorem using tropical intersection theory.
\end{abstract}

\section{Introduction}

One of the most celebrated results of tropical geometry certainly is Mikhalkin's Correspondence Theorem \cite{Mik05} which established a remarkable connection between enumerative algebraic geometry and tropical geometry. The result of the theorem is twofold: first, it showed that a difficult algebraic enumerative problem, namely that of counting plane curves of given degree and genus through an appropriate number of generic points, can also be formulated in the tropical world in a meaningful way and be solved combinatorially; and second, it proved that a solution of the tropical problem also solves the algebraic problem by showing that the multiplicities of the tropical curves are the numbers of algebraic curves tropicalizing to them. Both parts of the theorem motivated different kinds of questions in tropical geometry.

The first part, which basically laid the foundations for tropical enumerative geometry, led to the natural question of whether the methods used in the modern treatment of enumerative algebraic geometry could be used in tropical geometry as well. More precisely, the question was if it was possible to construct tropical moduli spaces of objects one wants to count and then count them using intersection-theoretical methods. It turned out that this in fact works out: the first tropical moduli space appearing in the literature was the tropical Grassmannian \cite{SS04}, which was used to construct moduli spaces of abstract and parametrized marked tropical curves in \cite{GKM09}. These spaces have enough structure to apply the intersection-theoretical methods of tropical geometry developed in \cite{AR10}.

The question motivated by the second part was, of course, that of how to generalize Mikhalkin's results to other enumerative problems. The main difficulty always is to find a strong relation between the algebraic objects one wants to count and their tropical analogues. The usual way to proceed is to assign to a given algebraic object a naturally constructed tropical object, generally called its tropicalization. Then by considering the fibers of the tropicalization, the numbers of algebraic and tropical objects having certain properties can be related. There are different ways of how to tropicalize a given algebraic object. One way, which is the one used in Mikhalkin's paper, is to consider degenerations. This approach also has been used in \cite{Shustin05}, \cite{NS06}, and \cite{Tyomkin12} to prove generalizations of Mikhalkin's result. Another way is to find a torus embedding of the algebraic moduli space of the objects in question such that its tropicalization (its non-Archimedean amoeba) is equal to the tropical moduli space of the corresponding tropical objects. That this is in fact possible sounds like a strong assumption, yet it has turned out to be possible in an increasing number of instances. Since enumerative numbers are defined by intersection numbers on the moduli spaces in both algebraic and tropical geometry, one needs information about how to relate tropicalizations of intersections to intersections of tropicalizations. This has been studied in detail in \cite{OP13} and \cite{OR11}.

The goal of this paper is to show that the previous knowledge about moduli spaces, intersection theory, and the tropicalization map is sufficient to reprove the correspondence theorem in a manner that can easily be generalized and applied to similar problems. What we need is that the algebraic moduli space is very affine, tropicalizes to the tropical moduli space after choosing a suitable torus embedding, and that the properties of the objects we want to count are encoded in finitely many torus morphisms to some appropriate tori, on the tropical as well as on the algebraic side. Furthermore, the tropical and algebraic torus morphisms should be compatible with tropicalization. These assumptions are motivated by the problem of counting marked parametrized curves, where the interesting enumerative numbers can be defined via the evaluation maps on the algebraic, as well as on the tropical side. We assume that the valued field we are working over is algebraically closed and has characteristic $0$ so that we are able to use Bertini arguments. The valuation on the other hand can be arbitrary. In particular, it does not need to be non-trivial.

The paper is organized as follows. In Section \ref{pSec:Preliminaries} we will review the basic results about tropical cycles and varieties, tropical intersection theory, and tropicalization needed in the subsequent sections.
Section \ref{cSec:Curves} is about marked parametrized curves in projective space. We will start by reviewing the necessary results about the tropical case, that is the definition of abstract and parametrized tropical curves and the construction of their moduli spaces. Afterwards we will examine the algebraic side. We will see that the correct objects to study are labeled parametrized marked rational curves. In the last part of Section \ref{cSec:Curves} we will show that their moduli spaces have canonical torus embeddings with respect to which they tropicalize onto the tropical moduli spaces. This yields a way to tropicalize algebraic curves. We will give an explicit description of the tropicalization of a curve and prove that tropicalization commutes with evaluating the marked points. 
In Section \ref{tSec:Relating Algebraic and Tropical Enumerative Geometry} we will abstract the situation of Section \ref{cSec:Curves} in order to prove a general correspondence theorem. We will assume that we want to count the number of elements in the preimages of an appropriate number of ``evaluation'' maps, which are given on both the algebraic and tropical side and commute with tropicalization.
Finally, in Section \ref{aSec:Applications} we will apply our results. Our first application will be our motivating example of rational curves in projective space of given degree where our correspondence theorem yields that the algebraic curves through an appropriate number of generic points tropicalize precisely to the tropical curves through the tropicalizations of the points. Afterwards, we will generalize this result to curves in arbitrary toric varieties, whose corresponding tropical curves will then have arbitrary tropical degree. We will also see that by fixing the labeled points in the boundary of the toric variety and its tropical toric counterpart we get a correspondence of relative invariants. Finally, we will show how to use the tropicalization of the incidence correspondence of lines and cubic surfaces in three-space to intersection-theoretically define multiplicities on tropical lines in smooth tropical cubic surfaces. We will then use our correspondence theorem to prove that these multiplicities respect the relative realizability of the lines and, in particular, sum up to $27$.

I would like to thank Andreas Gathmann for many helpful discussions and comments.

\section{Preliminaries}
\label{pSec:Preliminaries}

In this section we want to recall the definitions and some properties of the basic objects and operations in tropical geometry needed in the subsequent sections. Before we begin with a short review of tropical cycles and varieties, let us fix some general notation used throughout this paper. We will denote algebraic objects by Latin letters and tropical objects by fraktur letters. The algebraic objects will always be defined over an algebraically closed valued field $(K,\nu)$ of characteristic $0$, with valuation ring $(R,\mathfrak m)$, residue field $\kappa=R/\mathfrak m$, and value group $\VG=\nu(K^*)$. For a set $S$ we will denote by $\dPair S$ the set of pairs of distinct elements of $S$.

\subsection{Tropical Varieties}
\label{subSec:Tropical Varieties}

A \emph{tropical torus} is a finite dimensional real vector space $\ttorus$, together with a lattice $\nt\subset\ttorus$ such that the canonical map $(\nt)_{\mathds R}=\nt\otimes_{\mathds Z}\mathds R\rightarrow \ttorus$ is an isomorphism. Very often, one starts with the lattice $N$ and then sets $\ttorus=\nr$. This happens especially when relating algebraic and tropical tori: Given an algebraic torus $T$, the group $N$ of $1$-parameter subgroups (1-psgs) is a lattice and the tropical torus $\nr$ is called the tropicalization of $T$. For example, the tropical torus $\mathds R^n=(\mathds Z^n)_{\mathds R}$ is the tropicalization of the algebraic torus $\mathds G_m^n$, whose $1$-psgs correspond canonically to $\mathds Z^n$. For us, the most important tropical tori will be the tori $\tT r\coloneqq\mathds R^{r+1}/\mathds R \mathbf 1$ for $r\in\mathds N$, where $\mathbf 1$ is the vector with all entries equal to $1$, and $N_{\tT r}$ is the image of $\mathds Z^{r+1}$ under the canonical projection $\mathds R^{r+1}\rightarrow\tT r$. Its importance stems from the fact that it can be canonically identified with the tropicalization of the algebraic torus $\aT r\subset\PS r$ consisting of all points of $\PS r$ not contained in any coordinate hyperplane. Given two tropical tori $\ttorus$ and $\ttorus'$, a \emph{morphism} from $\ttorus$ to $\ttorus'$ denotes an $\mathds R$-linear map $\ttorus\rightarrow \ttorus'$ which is induced by a lattice homomorphism $\nt\rightarrow N_{\ttorus'}$.

Now we will review the definition of tropical cycles following \cite{AR10}. A \emph{weighted complex} in a tropical torus $\ttorus$ is a pure-dimensional rational polyhedral complex, together with a weight function assigning an integer to each of its inclusion-maximal cells. Here, rational means that each of its cells is a finite intersection of half-spaces $\{\mathfrak p\in\ttorus\mid \langle m,\mathfrak p\rangle\leq x\}$, where $m\in (\nt)^\vee=\Hom(\nt,\mathds Z)$ and $x\in\mathds R$. Suppose $\sigma$ is a rational polyhedron in $\ttorus$, and $\mathfrak p\in\sigma$. Then the subspace $U=\Lin(\sigma-\mathfrak p)$ is independent of the choice of $\mathfrak p$. Thus it makes sense to define $N_\sigma\coloneqq\nt\cap U$. It is a saturated sublattice of $\nt$ of rank $\dim(\sigma)$. If $\tau$ is a face of $\sigma$, then $N_\sigma/N_\tau$ is a lattice of rank $1$ and therefore has exactly two generators. The \emph{lattice normal vector} $u_{\sigma/\tau}$ is defined as that one of them pointing in the direction of $\sigma$. 
A \emph{tropical complex} $\mathfrak A$ is a weighted complex satisfying the so-called \emph{balancing condition}. This condition requires that whenever $\tau$ is a codimension $1$-cell of $\mathfrak A$, and $\sigma_1,\dots,\sigma_k$ are the maximal cells containing $\tau$, with weights $\omega_1,\dots,\omega_k$, respectively, then 
\begin{equation*}
\sum_{i=1}^k\omega_i u_{\sigma_i/\tau}= 0 \quad\text{in  }\nt/N_\tau.
\end{equation*}
The \emph{support} $|\mathfrak A|$ of $\mathfrak A$ is defined as the union of all its maximal cells with nonzero weight. 

\emph{Tropical cycles} are tropical complexes modulo refinements. Of course, refinements have to respect weights so that the support of a cycle is well-defined. For $d\in\mathds N$ let $Z_d(\ttorus)$ be the set of $d$-dimensional tropical cycles in $\ttorus$. There is a naturally defined addition of tropical cycles of the same dimension, making $Z_d(\ttorus)$ an abelian group. Given a cycle $\mathfrak A$ in $\ttorus$, the set of $d$-dimensional cycles $\mathfrak B\in Z_d(\ttorus)$ with $|\mathfrak B| \subseteq|\mathfrak A|$ is a subgroup of $Z_d(\ttorus)$, denoted by $Z_d(\mathfrak A)$. The graded group $\bigoplus_d Z_d(\mathfrak A)$ is denoted by $Z_*(\mathfrak A)$.

Let $\mathfrak A$ and $\mathfrak B$ be tropical cycles in tropical tori $\ttorus$ and $\ttorus'$, respectively. Then the \emph{morphisms} from $\mathfrak A$ to $\mathfrak B$ are all maps $|\mathfrak A|\rightarrow |\mathfrak B|$ which are induced by a morphism $\ttorus\rightarrow\ttorus'$ of tropical tori.

Finally, a \emph{tropical variety} is a tropical cycle with positive weights only.

\subsection{Tropical Intersection Theory}

Over the last few years tropical geometers have been working to develop a tropical intersection theory that parallels the algebraic intersection theory in many ways. As the name suggests, the tropical objects corresponding to the algebraic cycle groups of varieties are the groups of tropical cycles. Here, we want to give a very brief overview over the operations of tropical intersection theory needed in the rest of this paper.

\begin{itemize}
\item If $\ttorus$ is a tropical torus, there is an \emph{intersection product} on $Z_*(\ttorus)$, usually denoted by $``\cdot''$, making it into a commutative ring with unity $\ttorus$ \cite{AR10}. Similarly as in algebraic intersection theory, we define $Z^d(\ttorus)=Z_{\dim(\ttorus)-d}(\ttorus)$ and $Z^*(\ttorus)=\bigoplus_d Z^d(\ttorus)$. With this notation $Z^*(\ttorus)$ is even a graded ring.

\item Intersection products cannot only be defined for cycles in tropical tori, but more generally for cycles in \emph{smooth} tropical varieties, that is varieties locally isomorphic to the Bergman fan of a matroid. This has been shown in \cite{FR10} using similar methods as in \cite{AR10}, and in \cite{Sh13} using tropical modifications.

\item There is a concept of \emph{rational functions} on a tropical variety $\mathfrak A$. They are defined as continuous functions $\phi:|\mathfrak A|\rightarrow \mathds R$ such that there exists a polyhedral subdivision of $|\mathfrak A|$ on whose maximal cells $\phi$ is given by integer affine functions. 
Rational functions form a group with respect to addition which we denote by $\mathcal K(\mathfrak A)$. Every rational function $\phi\in\mathcal K(\mathfrak A)$ has an associated cycle $\phi\cdot\mathfrak A\in Z_{\dim(\mathfrak A)-1}(\mathfrak A)$ \cite{AR10}. The induced map 
\begin{equation*}
\mathcal K(\mathfrak A)\times Z_*(\mathfrak A)\rightarrow Z_*(\mathfrak A), (\phi,\mathfrak B)\mapsto \phi|_{|\mathfrak B|}\cdot\mathfrak B
\end{equation*}
is bilinear. Intersection with rational functions is also compatible with intersection products on smooth tropical varieties \cite{FR10,Sh13}.

\item If we consider sums of products of rational functions we obtain the concept of piecewise polynomials. These can be intersected with tropical cycles as well. They are the local building blocks for tropical cocycles, as introduced in \cite{F11}. On every tropical variety $\mathfrak A$, the tropical cocycles form a graded ring $C^*(\mathfrak A)$ and there is a bilinear map, analogous to the cap product in algebraic geometry,
\begin{equation*}
C^*(\mathfrak A)\times Z_*(\mathfrak A)\rightarrow Z_*(\mathfrak A)
\end{equation*}
which is graded in the canonical way, and compatible with intersections with rational functions. On smooth tropical varieties it is also compatible with intersection products, and in the special case of a tropical torus $\ttorus$ it is known that the map
\begin{equation*}
C^*(\ttorus)\rightarrow Z^*(\ttorus), \phi\mapsto \phi\cdot\ttorus
\end{equation*}
is an isomorphism of rings \cite[Thm 2.25]{F11}.

\item If $\mathfrak f:\mathfrak A\rightarrow\mathfrak B$ is a morphism, there is a \emph{push-forward} map $\mathfrak f_*: Z_*(\mathfrak A)\rightarrow Z_*(\mathfrak B)$ preserving the dimension of cycles. The push-forward is functorial \cite{GKM09,AR10}.

\item If $\mathfrak f$ is as above, there also is a \emph{pull-back} $\mathfrak f^*:C^*(\mathfrak B)\rightarrow C^*(\mathfrak A)$, which is a morphism of graded rings. In particular, if $\mathfrak B$ is a tropical torus, we can pull back tropical cycles. Pull-back and push-forward can be related in a tropical projection formula, that is whenever $\mathfrak C\in Z_*(\mathfrak A)$ and $\phi\in C^*(\mathfrak B)$ we have $\mathfrak f_*(\mathfrak f^*\phi\cdot\mathfrak C)=\phi\cdot\mathfrak f_*\mathfrak C$ \cite[Prop. 2.24]{F11}.

\item There is a notion of \emph{rational equivalence}. Let $\mathfrak A$ be a tropical variety. We denote by $\mathcal R(\mathfrak A)\subset\mathcal K(\mathfrak A)$ the set of all bounded rational functions on $\mathfrak A$. Two cycles $\mathfrak B,\mathfrak C\in Z_d(\mathfrak A)$ are said to be rationally equivalent if there is a tropical variety $\mathfrak D$ and a morphism $\mathfrak f:\mathfrak D\rightarrow \mathfrak A$ such that $\mathfrak B-\mathfrak C=\mathfrak f_*(\sum_i \phi_i\cdot\mathfrak E_i)$ for some $\phi_i\in\mathcal R(\mathfrak D)$ and $\mathfrak E_i\in Z_{d+1}(\mathfrak D)$. The canonical degree map $\deg:Z_0(\mathfrak A)\rightarrow \mathds Z$ is well-defined modulo rational equivalence \cite{AHR16}. 
\end{itemize}

\subsection{Tropicalization}
\label{pSubsec:Tropicalization}
Now we want to recall how to obtain a tropical variety from a subscheme of an algebraic torus. Let $T$ be an algebraic torus, let $N$ be the lattice of one-parameter subgroups of $T$, and let $M$ be the lattice of characters of $T$, which is canonically dual to $N$. Then there is a canonical isomorphism $T\cong\Spec (K[M])$. Any $K$-rational point $x\in T(K)$ induces a group homomorphism $M\rightarrow K^*$. Composing this with the valuation $\nu:K^*\rightarrow\mathds R$ we obtain an element in $\ttorus\coloneqq\nr\cong\Hom(M,\mathds R)$ called the tropicalization of $x$. By this procedure we obtain a map
\begin{equation*}
\trop:T(K)\rightarrow\ttorus.
\end{equation*}
It follows immediately from the construction that the tropicalization map factors through $\VG(\ttorus)\coloneqq\VG\otimes_{\mathds Z} N$. 

Now let $X$ be a pure-dimensional subscheme of $T$ and $\mathfrak t\in \VG(\ttorus)$. Then there is a $t\in T$ with $\trop(t)=\mathfrak t$. The \emph{initial degeneration} $\ini_{\mathfrak t}X$ is defined as the special fiber of the scheme-theoretic closure of $t^{-1}X$ in $\mathcal T=\Spec (R[M])$. It is well-defined up to translation by a point in $T_\kappa=\Spec(\kappa[M])$. If $\mathfrak t$ is not in $\VG(\ttorus)$ we can use the same definition after a suitable extension of valued fields. Of course, this depends on the choice of the field extensions, but it can be shown that for every two of those extensions $L$ and $L'$ there is a common field extension $\overline L$ such that the initial degenerations defined over $L$ and $L'$ become translates of each other after passing to $\overline L$ \cite[Section 5]{Gub12}. We remark that it is also possible to define initial degenerations at general $\mathfrak t\in\ttorus$ using the tilted group rings $R[M]^{\mathfrak t}$. With this definition they are well-defined without the choices of a field extension and a $t\in T$ \cite{OP13}.

The \emph{tropicalization} $\Trop(X)$ of $X$ can be defined as a tropical cycle in $\ttorus$ whose underlying set $|\Trop(X)|$ consists of the $\mathfrak t\in\ttorus$ with $\ini_{\mathfrak t}X\neq\emptyset$. This is indeed a purely $\dim(X)$-dimensional polyhedral complex \cite{BG84}, and if the valuation is non-trivial it is the closure of $\trop(X(K))$ \cite{SS04, Drai08,P09,OP13}. For any point $\mathfrak t\in|\Trop(X)|$, the sum of the multiplicities of the components of $\ini_{\mathfrak t} X$ is well-defined and called the \emph{tropical multiplicity} of $X$ at $\mathfrak t$. The multiplicity is constant on the relative interiors of the maximal cells of any polyhedral structure on $|\Trop(X)|$ and satisfies the balancing condition \cite{Speyer-diss,ST08,Gub12}. Thus, we get a well-defined tropical cycle $\Trop(X)$.

The fiber of $0$ with respect to the tropicalization map $\trop$ can easily be described algebraically. It is exactly the image of the inclusion $\mathcal T(R)\hookrightarrow T(K)$, which maps a morphism $\Spec R\rightarrow \mathcal T$ to the induced morphism $\Spec K\rightarrow T$ on the generic fibers. For this reason we call points in $T(K)$ tropicalizing to $0$ the $R$-\emph{integral points} of $T$. We can use this identification to give $\trop^{-1}\{0\}$ a topology. Namely, we give it the initial topology with respect to the reduction map $\mathcal T(R)\rightarrow T_\kappa(\kappa)$, where we consider $T_\kappa(\kappa)$ with the Zariski topology. This makes $\trop^{-1}\{0\}$ a topological group. Note that if the valuation is trivial, we have $T(K)=\mathcal T(R)=T_\kappa(\kappa)$ and the topology just defined is the ordinary Zariski topology on $T(K)$. 

If $\mathfrak t\in\VG(\ttorus)$ is general, then the fiber $\trop^{-1}\{\mathfrak t\}$ is a $\trop^{-1}\{0\}$-torsor and therefore has a well-defined topology induced by the topology on $\trop^{-1}\{0\}$. It could also be defined analogously to the topology of the fiber over $0$ using tilted group rings.

\section{Parametrized Rational Curves and their Moduli Spaces}
\label{cSec:Curves}

In this section we will review the definitions of tropical and algebraic rational curves of certain degree. We will show that the algebraic moduli space of curves can be embedded into a torus in such a way that it tropicalizes to the tropical moduli space of curves and that this respects the evaluation maps. Furthermore, we show that the algebraic evaluation maps are induced by morphisms of tori. This will be the motivation for the abstraction of the problem of relating algebraic and tropical enumerative numbers considered in the next section.

\subsection{Tropical Curves}
\label{cSubSec:Tropical Curves}

An \emph{abstract (rational) tropical curve} is a connected metric graph $\mathfrak C=(V,E,\ell)$ of genus $0$, $V$ being the set of vertices, $E$ the set of edges, and $\ell:E\rightarrow \mathds R_{>0}\cup\{\infty\}$ a length function. The edges adjacent to a $1$-valent vertex are called \emph{legs}, and we require that the legs are exactly those edges with infinite lengths. All other edges are called \emph{bounded edges}. The $1$-valent vertices should be imagined as the points at infinity of the unique leg incident to them and  will be called \emph{feet}. All other vertices will be referred to as \emph{inner vertices}. The interior part $\mathfrak C^\circ$ of $\mathfrak C$ will denote the set of points of $\mathfrak C$ not equal to a foot. As every leg should only have one foot, we exclude the pathological case of graphs consisting of exactly $2$ vertices. In that way, every leg $e\in E$ is incident to a unique inner vertex $v_e$. If an abstract tropical curve has a $2$-valent vertex, we can replace it and the two edges incident to it by one edge of adequate length. The resulting curve should obviously be considered isomorphic to the original one. For this reason, we restrict our attention to abstract tropical curves without $2$-valent vertices. The advantage of this is that in this way the category of abstract tropical curves, whose morphisms are morphisms of the underlying graphs respecting the lengths, induces the ``correct'' notion of isomorphisms.

If $I$ is some finite index set, we define an \emph{$I$-marked tropical curve} to be an abstract tropical curve together with an injection $I\rightarrow E,i\rightarrow e_i$ assigning a leg to every index in $I$. Of course, morphisms of marked curves have to respect the markings. We denote by $\tMno I$ the space of all $I$-marked  tropical curves with exactly $|I|$ legs, modulo isomorphisms.

The set $\tMno I$ can be identified with a tropical variety in the following way. For every tropical curve $(\mathfrak C,(e_i)_{i\in I})\in \tMno I$ and pair of distinct indices $i,j\in I$ there is a well-defined distance $\dist_{\mathfrak C}(i,j)$ between the vertices $v_{e_i}$ and $v_{e_j}$. We use these distances to define the tropical Plücker embedding which maps $\mathfrak C$ to the point 
\begin{equation*}
\tpl(\mathfrak C)\coloneqq\left(-\dist_{\mathfrak C}(i,j)/2\right)_{(i,j)}\in \mathds R^\dPair I/\Image\Phi,
\end{equation*}
where $\Phi$ is the morphism
\begin{equation*}
\Phi:\mathds R^I\rightarrow \mathds R^\dPair I, (x_i)_i\mapsto (x_i+x_j)_{(i,j)}.
\end{equation*}

It follows from the results of \cite[Thm. 4.2]{SS04} and \cite[Thm. 3.7]{GKM09} that $\tpl$ is an embedding and that its image is a pure-dimensional fan satisfying the balancing condition after giving weight $1$ to all inclusion-maximal cones. Note that both, Speyer and Sturmfels, and Gathmann, Kerber, and Markwig, consider slightly different embeddings. Namely, they consider unordered pairs instead of ordered ones, and they leave out the factor $-1/2$ in front of the distances. Our reason to work with ordered pairs is that they are more natural on the algebraic side, as we will see later.  The factor $-1/2$ is included because the lattice  in $\mathds R^\dPair I/\Image\Phi$ used in \cite{GKM09} is not the canonical one, which is induced by $\mathds Z^\dPair I\subseteq\mathds R^\dPair I$. It follows from \cite[Example 7.2]{FR10}, where it is shown that $\tMno I$ is a tropical matroid variety, that they instead use exactly $2$ times the canonical lattice.

Now we want to define parametrized $I$-marked (rational) tropical curves. To do this, we first need the notion of degree. Let $\ttorus$ be a tropical torus, and $J$ a finite set (disjoint from $I$). Then every map $\Delta:J\rightarrow \nt$ such that $\sum_{j\in J}\Delta(j)=0$ is called a \emph{degree} (in $\ttorus$ with index set $J$). A \emph{(labeled) parametrized $I$-marked (rational) tropical curve in $\ttorus$ of degree $\Delta$} is a triple $(\mathfrak C,(e_\lambda),\tpara)$ consisting of an $L_0\coloneqq I\cup J$-marked curve $(\mathfrak C,(e_\lambda))$ and a map $\mathfrak C^\circ\rightarrow \ttorus$ which is linear on the edges (that is affine linear after identifying an edge $e$ with $[0,\ell(e)]$ if it is bounded and with $[0,\infty)$ in case it is a leg) and satisfies the additional requirements stated in the following. For each inner vertex $v$ incident to an edge $e$ there is a well-defined direction $\dir(v,e)$, which is the derivative of $\mathfrak f|_e$ when choosing the parametrization of $e$ starting at $v$ . In case $e$ is a leg, we also denote $\dir(v,e)$ by $\dir(e)$. We require that all these direction vectors are in $\nt$, and, furthermore, that at every inner vertex $v$ the balancing condition is fulfilled, that is
\begin{equation*}
\sum_{e:v\in e} \dir(v,e)=0.
\end{equation*}
Finally, we require that $\dir(e_j)=\Delta(j)$ for all $j\in J$, and $\dir(e_i)=0$ for all $i\in I$.
Morphisms of parametrized curves should, of course, respect the maps into $\ttorus$. We denote by $\tMlno{I}{\ttorus}{\Delta}$ the set of all parametrized $I$-marked tropical curves in $\ttorus$ of degree $\Delta$ having exactly $|L_0|$ legs, modulo isomorphisms.

For every $i\in I$ there is an evaluation map $\tev i:\tMlno{I}{\ttorus}{\Delta}\rightarrow \ttorus$ which assigns to a parametrized tropical curve the unique point in $\ttorus$ to which $e_i$ is mapped. 

We can identify $\tMlno{I}{\ttorus}{\Delta}$ with a tropical variety in a similar way as $\tMno I$. To do so, fix $i_0\in I$. Then by \cite[Prop. 4.7]{GKM09} the map
\begin{equation*}
\tpl\times\tev {i_0}:\tMlno{I}{\ttorus}{\Delta}\rightarrow \mathds R^\dPair {L_0}/\Image\Phi\times\ttorus,
\end{equation*}
is an embedding with image $\tMno {L_0}\times\ttorus$.

The most important situation for us is when $\ttorus=\tT r$. In this case, we obtain a degree $\Delta_d$ for every $d\in\mathds N$ by setting $J=\{0,\dots,r\}\times \{1,\dots, d\}$ and $\Delta((i,j))=s_i$, where $s_i$ denotes the image of the $i$-th standard basis vector of $\mathds R^{r+1}$ in $\tT r$. By abuse of notation, we usually just write $d$ instead of $\Delta_d$. The reason why this case is of special importance is that $\tT r$ is the tropicalization of $\aT r$, so that the tropicalization of an algebraic curve of degree $d$ in $\PS r$ which intersects the coordinate hyperplanes generically should intuitively be some tropical curve of degree $d$ in $\tT r$.

\subsection{Algebraic Curves}
\label{cSubSec:Algebraic Cuves}

An \emph{abstract (rational) (algebraic) curve} is a curve isomorphic to $\PS 1$. If $I$ is a finite index set, then an \emph{$I$-marked abstract curve} is is constituted by an abstract curve $C$ together with an injective map $I\rightarrow C$. 
In the same vein as in the tropical case we continue with the definition of parametrized $I$-marked (rational) curves in $\PS r$ of certain degree. In the algebraic setting, as degree we allow elements in $A_1(\PS r)\cong \mathds Z$. For $d\in\mathds Z$, a \emph{parametrized $I$-marked (rational) curve in $\PS r$ of degree $d$} consists of an $I$-marked abstract curve $C$ together with a morphism $\apara:C\rightarrow\PS r$ such that $\apara_*[C]=d$. It can be shown that there is a moduli space $\Mno I$ whose points correspond bijectively to isomorphism classes of theses curves. It is easily seen that for $i\in I$ there are natural evaluation maps
\begin{equation*}
\ev i:\Mno I\rightarrow \PS r
\end{equation*}
assigning to a parametrized $I$-marked curve the image of the $i$-th marked point, and again it can be shown that these maps are in fact morphisms of algebraic varieties. With these maps, counting curves through given distinct points $p_1,\dots,p_n\in\PS r$ is, morally speaking, the same as determining the number of elements of 
\begin{equation}
\label{cEqu:Preimage of Points}
\ev 1^{-1}\{p_1\}\cap\dots\cap \ev n^{-1}\{p_n\}\subseteq\Mn,
\end{equation}
given, of course, that there are only finitely many.

We would like to count curves using tropical methods, that is we would like to replace the evaluation maps $\ev i$ in Expression \ref{cEqu:Preimage of Points} by the tropical evaluation maps $\tev i$,  the moduli space $\Mn$ by the tropical moduli space $\tMln$, and the points by appropriate points in $\tT r$. Of course, it is by no means clear that the tropical enumerative numbers arising this way are equal to the algebraic ones. On the contrary, it is a deep result first proven in \cite{Mik05}, at least for plane curves. We want to reprove this statement by relating the algebraic and tropical enumerative problems via the tropicalization map. For this to work, we need to find an embedding of $\Mn$ into some torus such that $\tMln$ is its tropicalization. However, this cannot be expected to work for two reasons: first, the curves described by $\Mn$ map into projective $r$-space $\PS r$, whereas the tropical curves map into $\tT r$, which corresponds to $\aT r$. Therefore, it is unclear how to tropicalize the points in the boundary of $\PS r$. Second, tropical curves comprise labels of additional ends, whereas algebraic curves lack this structure. These shortcomings can be remedied in a single step as we will see in the following.

A general curve $\mathcal C =(C,(p_i)_{i\in I}, \apara)$ in $\Mno I$ is expected to intersect the coordinate hyperplanes $H_i=\{x_i=0\}$ of $\PS r$ in exactly $d$ unmarked points each, and no point of the curve is supposed to lie in more than one of those hyperplanes. Hence there are $d(r+1)$ distinct points $p_{ij}\in C$, for $(i,j)\in J\coloneqq\{0,\dots, r\}\times\{1,\dots, d\}$ such that $\apara^{-1}H_i$ consists exactly of the points $p_{i1}\dots p_{id}$. The images $\apara(p_{ij})$ are exactly those points where the curve $\apara(C)$ meets the boundary divisors of the toric variety $\PS r$, so the points $p_{ij}$ should correspond to the feet of the tropicalization of the intersection of $\apara(C)$ with the torus $\aT r\subseteq \PS r$.  This motivates the definition of a \emph{labeled parametrized $I$-marked (rational) curve in $\PS r$ of degree $d$} as an $L_0=I\cup J$-marked parametrized curve $\mathcal C=(C,(p_\lambda)_{\lambda\in L_0},\apara)$ in $\PS r$ of degree $d$ such that $\apara^{-1}H_i=\{p_{ij}\mid 1\leq j\leq d\}$. We denote by $\Mlno I\subseteq \Mno {L_0}$ the set of all isomorphism classes of these labeled curves. Note that for $i\in I$ the evaluation maps $\ev i:\Mno {L_0}\rightarrow \PS r$ induce evaluation maps $\ev i:\Mlno I\rightarrow \aT r$.

\subsection{Tropicalizing $\Mlno I$}
\label{cSubSec:Tropicalizing}

Now we want to show that $\Mlno I$ is the right algebraic object to study, that is that it tropicalizes onto $\tMlI I$ in a way respecting the evaluation maps. We first need to find a suitable embedding of $\Mlno I$ into an algebraic torus. To make the transition from algebraic to tropical geometry as clear as possible we use the same notation as in \ref{cSubSec:Tropical Curves}, that is we choose an element $i_0\in I$ and write $L_0= I\cup J$ (where we still denote $J=\{0,\dots, r\}\times\{1,\dots, d\}$). Let  $\mathcal C=(C, (p_\lambda)_{\lambda\in L_0}, \apara)\in\Mlno I$. As $C$ is only well defined up to isomorphisms we can assume $C=\PS 1$, in which case we can write $p_\lambda=(b_\lambda:a_\lambda)$ in homogeneous coordinates for appropriate $a_\lambda, b_\lambda\in K$. Since the curve has degree $d$, the morphism $\apara$ is given by $r+1$ homogeneous polynomials $f_0,\dots, f_r$ of degree $d$. The fact that $\apara^{-1}H_i= \{p_{i1},\dots,p_{id}\}$ implies that $f_i$ is of the form $f_i= c_i\prod_{1\leq j\leq d} (a_{ij}x-b_{ij}y)$ for some $c_i\in K^*$. Because $p_{i_0}$ is not mapped to any $H_i$ it follows that $\apara$ is completely determined by the $p_{ij}$ and $\ev {i_0}(\mathcal C)$. Let $U\subset(\PS1)^{L_0}$ be the complement of the big diagonal, that is the set of $L_0$-tuples of distinct points in $\PS1$. By what we just saw, there is a surjection $\phi:U\times \aT r\rightarrow\Mlno {L_0}$ mapping an element $((p_\lambda)_{\lambda\in L_0},P)$ to the curve $(\PS1, (p_\lambda)_{\lambda\in L_0},\apara)$, where $\apara$ is the unique morphism mapping $p_{ij}$ into $H_i$ and $p_{i_0}$ to $P$. If we consider the action of $\GL(2)$ on $U\times \aT r$ induced by the diagonal action on $(\PS1)^{L_0}$, it is immediate that this surjection induces an identification of $\Mlno I$ with $U\times \aT r/\GL(2)$. 

By the Gelfand-MacPherson correspondence \cite{GM82} the quotient $U/\GL(2)$ is in bijection with the subset $\G^0(2,L_0)$ of the Grassmannian $\G(2,L_0)$ consisting of all lines that do not pass the intersection of two coordinate  hypersurfaces, modulo the canonical action of the torus $\mathds G_m^{L_0}$. Let us quickly review the bijection: An orbit $x$ in $U/\GL(2)$ is represented by a tuple $(p_\lambda)_{\lambda\in L_0}$. For each of the $p_\lambda$ we take a representation $(b_\lambda: a_\lambda)$ in homogeneous coordinates. Writing these coordinates into the columns of a matrix we obtain a $2\times L_0$-matrix of rank $2$. This matrix gives us coordinates for a point in $\G^0(2,L_0)$. The Gelfand-MacPherson correspondence takes $x$ to the orbit of this point.

Now we consider the Plücker embedding of $\G(2,L_0)$ into $\mathds P K^ {\dPair {L_0}}$. Notice that here we do not use the usual Plücker embedding, which involves only one coordinate per minor, but we get each minor twice with opposite signs. In this way we do not need to choose any signs for the minors. The $\mathds G_m^{L_0}$-action on $\G(2,L_0)$ can be extended to $\mathds P K^{\dPair{L_0}}$ by setting $(x_\lambda) .(p_{(\lambda,\mu)})=(x_\lambda x_\mu p_{(\lambda,\mu)})$. The torus $\mathds T K^ {\dPair {L_0}}$ is invariant with respect to this action, and the induced action on it is induced by the morphism of tori given by $(x_\lambda)_\lambda\mapsto (x_\lambda x_\mu)_{(\lambda,\mu)}$. The quotient $\mathds T K^ {\dPair {L_0}}/\mathds G_m^{L_0}$ is an algebraic torus again, and since the action on $\mathds T K^ {\dPair {L_0}}$ is induced by the natural $\mathds G_m^{L_0}$-action on $\mathds G_m^{\dPair{L_0}}$ it is equal to the quotient $\mathds G_m^{\dPair{L_0}}/\mathds G_m^{L_0}$. Thus we have obtained an embedding of $\Mlno I$ into the torus $\mathds G_m^{\dPair{L_0}}/\mathds G_m^{L_0}\times\aT r$.

This embedding can be used to compute a tropicalization of $\Mlno I$. Actually, the result is well-known: $\Mlno I$ is equal to the image in $\mathds T K^ {\dPair{L_0}}/\mathds G_m^{L_0}$ of the Plücker-embedding of $\G^0(2,L_0)$ into $\mathds T K^{\dPair{L_0}}$, times $\aT r$. Hence its tropicalization is equal to $\tMno {L_0}\times \tT r$ by \cite[Thm. 3.4]{SS04}, which in term is equal to $\tMlI I$ by \cite[Prop. 4.7]{GKM09}. We see that the tropicalization of $\Mlno I$ is exactly its tropical analogue. However, to really be sure that we have chosen the correct torus embedding we need to show that tropicalization commutes with evaluation maps. We begin by constructing a natural candidate for the tropicalization of a given curve $\mathcal C=(C, (p_\lambda)_{\lambda\in L_0}, \apara)$ in $\Mlno I$, which will be seen to be the tropical curve corresponding to the tropicalized Plücker coordinates of $\mathcal C$. For this we assume that we are given a representative of $\mathcal C$ in \emph{standard form}, that is that $C=\PS 1$ and $p_{i_0}=(0:1)$. Our construction is a generalization of \cite[Constr. 2.2.20]{DennisDiss} to general valued fields. It works for parametrized curves in arbitrary toric varieties, as explained in Section \ref{aSec:Applications}, yet we will restrict ourselves to curves in projective space here to reduce the amount of notation.

\begin{construction}
\label{cConstr:Corresponding tropical curve}
Let $\mathcal C=(\PS 1, (p_\lambda)_{\lambda\in L_0}, \apara)$ be a labeled parametrized $I$-marked curve in $\PS r$ of degree $d$ in standard form. We will construct a parametrized $I$-marked tropical curve in $\tT r$ of degree $d$ as a candidate for the tropicalization of $\mathcal C$. Since $\mathcal C$ is in standard form we have $p_{i_0}=(0:1)$. Because all marked (resp. labeled) points are distinct, there are unique $a_\lambda\in K^*$ such that $p_\lambda=(1:a_\lambda)$ for $\lambda\neq i_0$. Furthermore, suppose that $\ev {i_0}(\mathcal C)=\apara(p_{i_0})=(c_0:\dots:c_r)$.

For each nonempty $A\in\mathcal P(L)\setminus\{\emptyset\}$ in the power set of $L\coloneqq L_0\setminus \{i_0\}$ we define $\nu(A)= \min\{\nu(a_\mu-a_\lambda)\mid \lambda,\mu\in A\}\in\mathds R\cup\{\infty\}$. We define a partial order on $\mathcal P(L)\setminus\{\emptyset\}$ by letting $A\preceq B$ if and only if $A\subseteq B$ and $\nu(A)= \nu(B)$. Let $V\subseteq\mathcal P(L)\setminus\{\emptyset\}$ be the set of all maximal elements with respect to this order. Note that all singletons belong to $V$. Ordering $V$ by inclusion, we obtain another partially ordered set. We make it into the underlying set of a graph $G$ by connecting two elements of $V$ if and only if one covers the other. As $L$ is greater than every other element of $V$, there is a path in $G$ between any vertex $v\in V$ and $L$. This shows that $G$ is connected. Even more, $G$ is a tree: If $A, B$ are two subsets of $L$ with nonempty intersection, then it follows from the valuation properties that $\nu(A\cup B)=\min\{ \nu(A),\nu(B)\}$. Thus we either have $A\preceq A\cup B$ or $B\preceq A\cup B$. If both $A$ and $B$ are in $V$ then it follows from the definitions that either $A\subseteq B$ or $B\subseteq A$. In particular, every vertex $v\neq L$ is covered by precisely one element in $V$. This shows that $G$ has exactly $|V\setminus\{L\}|=|V|-1$ edges, which implies that it is a tree. 

We make $G$ into an abstract tropical curve $\mathfrak C$ by giving an edge $e=\{v\subset w\}$ the length $\ell(e)=\nu(v)-\nu(w)$ and attaching an extra foot $\{i_0\}$ along a leg to $L$. This curve is $L_0$-marked in a canonical way. We use $\mathfrak C$ to construct a labeled parametrized $I$-marked tropical curve in $\tT r$ of degree $d$. All there is left to do is to define the morphism $\tpara:\mathfrak C^\circ\rightarrow \tT r$. To do so, let us introduce some more notation. Let $\pi:J\rightarrow \{0,\dots,r\}$ be the projection onto the first coordinate. Furthermore, denote by $s_i$ the image of the $i$-th standard basis vector of $\mathds R^{r+1}$ in $\tT r$. For every subset $A\subset L_0$ we define 
\begin{equation}
\label{cEqu:Definition of s_a}
s_A\coloneqq \sum_{i=0}^r |\pi^{-1}\{i\}\cap A|s_i.
\end{equation}
Let $v\in V$ be an inner vertex, and let $L=v_0, v_1, \dots, v_k=v$ be the vertices passed by the unique path from $L$ to $v$. We define 
\begin{equation*}
q_v\coloneqq\sum_{i=0}^r\nu(c_i)s_i+\sum_{i=1}^k \ell(\{v_i,v_{i-1}\})s_{v_i}
\end{equation*}
and use this to define $\tpara$ as the map $\mathfrak C^\circ\rightarrow \tT r$ which sends an inner vertex $v$ to $q_v$, interpolates linearly on the bounded edges, and has direction $s_v$ on a leg with foot $v$. To ensure that this really defines a tropical curve we need to check the balancing condition. It is immediate that for any edge $e=\{v\subset w\}$ in $\mathfrak C$ we have $\dir(w,e)=s_v$. If $v$ is an inner vertex, and $w=v_0,\dots,v_k$ are the vertices adjacent to $v$, then $v$ is the disjoint union of $v_1,\dots, v_k$. This yields

\begin{equation*}
\sum_{i=0}^k\dir(v,\{v,v_i\})= -s_v+\sum_{i=1}^k s_{v_i}=-s_v+s_{\bigcup_i v_i}=0,
\end{equation*}
that is the balancing condition is satisfied at $v$. For $v=L$ there is no $v_0$, but instead we can use $s_L=0$ to see that the balancing condition is satisfied here as well. As our curve has degree $d$ by construction, we really obtained a labeled parametrized $I$-marked tropical curve $\mathfrak C$ of degree $d$. We denote it by $\mathfrak C_{\mathcal C}\coloneqq\mathfrak C$ and call it the \emph{corresponding tropical curve to $\mathcal C$}.
\end{construction}

\begin{example}
\label{cEx:Corresponding tropical curve}
Let $K=\overline{\mathds C((t))}=\bigcup_{n\in\mathds Z}\mathds C((t^{\frac 1n}))$ be the field of Puiseux series. Furthermore, let $I=\{1,2\}$, $i_0=1$, and $d=r=2$. We want to compute the corresponding tropical curve to the curve $\mathcal C$ in standard form which is given by
\begin{alignat*}{4}
c_0&=1+t 	&\quad\qquad\quad	a_{01}&=	t^{-2}+1		&\qquad	a_{11}&=t^{-1}	&\qquad	a_{21}&=t^{-2} \\
c_1&=2t^{-2}+3t		&	a_{02}&=2			&	a_{12}&=2+t+4t^3	&	a_{22}&=2+t\\
c_2&=t^{-1}+t		&	a_{2\hphantom 1}&=2+t+4t^3-t^4		
\end{alignat*}
in the notation of the preceding construction. We can easily determine the elements of $V$ using the following three observations. First, the real numbers which are of the form $\nu(v)$ for some $v\in V$ are exactly those which are equal to $\nu(a_\lambda-a_\mu)$ for some $\lambda,\mu\in L$. Second, whenever $v\in V$ and $\lambda\in v$ we can reconstruct $v$ from $\lambda$ and $\nu(v)$ since $v=\{\mu\in L\mid \nu(a_\lambda-a_\mu)\geq\nu(v)\}$. Finally, for every $r\in\mathds R\cup\{\infty\}$ and $\lambda\in L$ we have $\{\mu\in L\mid \nu(a_\lambda-a_\mu)\geq r\}\in V$. 

We begin by computing the set $R\coloneqq\{\nu(a_\lambda-a_\mu)\mid \lambda,\mu\in L\}$. In our example case it is equal to $\{-2,-1,0,1,3,4,\infty\}$. We will obtain all elements of $V$ if we determine the equivalence classes under the equivalence relations
\begin{equation*}
\lambda\sim_r\mu:\Leftrightarrow \nu(a_\lambda-a_\mu)\geq r
\end{equation*}
for $r\in R$. Of course, it is advisable not to treat the computations of $L/\hspace{-3pt}\sim_r$ for different $r$ as independent problems. Sorting the elements of $R$ in an ascending order and computing the valuations of the vertices along the way will speed up the process considerably.
In our example, we start with $r=-2$ and, of course, obtain the vertex $L$ which has valuation $-2$. Then we continue with $r=-1$, where we obtain the two vertices $\{(0,1),(2,1)\}$ and $\{(1,1),(0,2),(1,2),(2,2),2\}$ with valuations $0$ and $-1$, respectively. That $\{(0,1),(2,1)\}$ has valuation $0$ tells us that it will also appear as equivalence class for $r=0$, the others being $\{(1,1)\}$ and $\{(0,2),(1,2),(2,2),2\}$ with valuations $\infty$ and $1$, respectively. Continuing like this, we see that we get the abstract tropical curve depicted in Figure \ref{cFig:curve} on the left. The lengths of the bounded edges can be computed directly from the valuations of their vertices. The directions of the edges can be read off from their vertices themselves. The image of $\mathfrak C_{\mathcal C}$ is also depicted in Figure \ref{cFig:curve}. Interestingly, it looks like a union of two tropical lines, even though it is the tropicalization of an irreducible curve. In fact, the parametrized tropical curve $\mathfrak C_{\mathcal C}$ can be considered as two $1$-marked tropical curves of degree $1$, glued together along a contracted edge.

\begin{figure}
\begin{subfigure}{.5\textwidth}
\flushleft
\begin{tikzpicture} [
dot/.style = {color=black,circle, fill=black, inner sep=0, minimum size=0.1cm},
nonexistent/.style ={inner sep=0, minimum size=0, fill=none},
label position=left]

\begin{scope}[ yscale=-0.7,xscale=1.1, font=\footnotesize]
\begin{scope}[every node/.style=dot,xshift=2.6cm, font=\tiny]
\draw (0,2) -- (0,3) node (1) [label={$\{1\}$}]{};
\node (L) at (0,2) [label=$L$] {}
	child [level distance =2cm, sibling distance= 2.1cm]{
	node (01 21 ) [label={$\{(0,1),(2,1)\}$}]{} 
		[level distance =5 cm, sibling distance= 0.7cm]
		child {node 	(01) [label={[label distance=-0.2cm]above:{$\{(0,1)\}$}}] {} }
		child {node 	(21) [label={[label distance=-0.2cm]above:{$\{(2,1)\}$}}] {} } 
		edge from parent
			node[nonexistent, left=0.07]{$2$}
	}
	child [level distance=1cm, sibling distance=1.4cm] {
	node (11 02 etc) [label=right:{$\{(1,1),(0,2),(1,2),(2,2),2\}$}]{}
		child [level distance=6cm, sibling distance =1.4cm]{
		node (11) [label={[label distance=-0.2cm]above:{$\{(1,1)\}$}}] {}
		}
		child [red,level distance =2cm,sibling distance=1.4cm]{
		node (02 12 etc) [label=right:{$\{(0,2),(1,2),(2,2),2\}$}] {}
			child [blue,level distance=4cm, sibling distance= 1.4cm] {
			node (02) [label={[label distance=-0.2cm]above:{$\{(0,2)\}$}}] {}
			}
			child [blue,level distance=2 cm, sibling distance=0.7cm] {
			node (12 22 etc) [label=right:{$\{(1,2),(2,2),2\}$}] {}
				child [level distance= 2cm, sibling distance= 0.7cm] {
				node (22) [label={[label distance=-0.2cm]above:{$\{(2,2)\}$}}] {}
				} 
				child [level distance= 1cm, sibling distance= 1.05cm] {
				node (12 2) [label=right:{$\{(1,2),2\}$}] {}
				[level distance=1cm]
					child [sibling distance= 0.35cm]{
					node (12) [label={[label distance=-0.2cm]above:{$\{(1,2)\}$}}] {}
					}
					child [sibling distance=1.05cm]{
					node (2) [label={[label distance=-0.05cm]above:{$\{2\}$}}] {}
					}
					edge from parent
						node[nonexistent, pos=0.65, left=0.05]{$1$}
				}
				edge from parent
					node[nonexistent, left=0.05]{$2$}
			}
			edge from parent
				node[nonexistent, left=0.05]{$2$}
		}
		edge from parent
			node[nonexistent, left=0.05,pos=0.6]{$1$}
	};
\end{scope}

\foreach \x in {-2,-1,...,4} {
 \draw (0.1,-\x) -- (-0.1,-\x) node [anchor=east] {$\x$}; 
 }
\draw (0,2.2) -- (0,-4.2);
\draw[dotted] (0,-4.2)-- (0, -4.8);
\draw (0,-4.8) -- (0,-5) node  [dot,label=$\infty$] {};
\draw[dotted] (0,2.2)-- (0, 3);
\node at (0,3.2) {Valuation};
\end{scope}
\end{tikzpicture}
\end{subfigure}
\begin{subfigure}{.49\textwidth}
\flushright
\begin{tikzpicture}[
dot/.style = {anchor =base,circle, fill=black, inner sep=0, minimum size=0.1cm},
nonexistent/.style ={anchor=base,inner sep=0, minimum size=0, fill=none},
label position=left]

\begin{scope}[xshift=10cm, yshift= 0cm,scale=0.4, font=\tiny,left, midway]
\clip (-6,-3) rectangle (5,5);

	\node (mark 1) at (-1,0) [dot,label={[label distance = -0.1cm]above left:$1$}] {};
	\draw (mark 1) -- (-3,0) node [nonexistent] (vert 1) {}; 
	\draw (vert 1) --node[pos=0.555,right=-0.03cm] {$(0,1)$}   +(-5,-5); 
	\draw (vert 1) --node[pos=0.91,left=-0.1cm]{$(2,1)$} +(0,5);
	
	\draw (mark 1) -- ++(1,0) node (contracted edge)[nonexistent]{} --node[pos=0.9,above=-0.1cm]{$(1,1)$} +(5,0);
	\draw [blue](contracted edge) +(-5,-5) --node[pos=0.44,black,right=-0.03cm]{$(0,2)$} +(0,0) -- ++(2,2) node [nonexistent] (vert 2) {};
	
	\node [dot,red] at (0,0) {};

	\draw [blue](vert 2) --node[black,right=-0.1cm]{$(2,2)$} +(0,5);
	\draw [blue](vert 2) -- ++(1,0) node (mark 2) [dot,label={[black,label distance =-0.1cm]above left:$2$}] {} --node [pos=0.305,black,above=-0.1cm] {$(1,2)$} +(5,0) ;	

\draw [->] (1,-3) -- (1,5);
\draw [->] (-6,1) -- (5,1);

\foreach \x in {-5,-4,...,4} 
\draw (\x,0.92)--(\x,1.08);

\foreach \x in {-2,-1,...,4}
\draw (0.92,\x)--(1.08,\x);
\end{scope}
\end{tikzpicture}
\end{subfigure}
\caption{The abstract tropical curve of $\mathfrak C_{\mathcal C}$ and its image in $\{0\}\times\mathds R^2\cong\tT 3$}
\label{cFig:curve}
\end{figure}
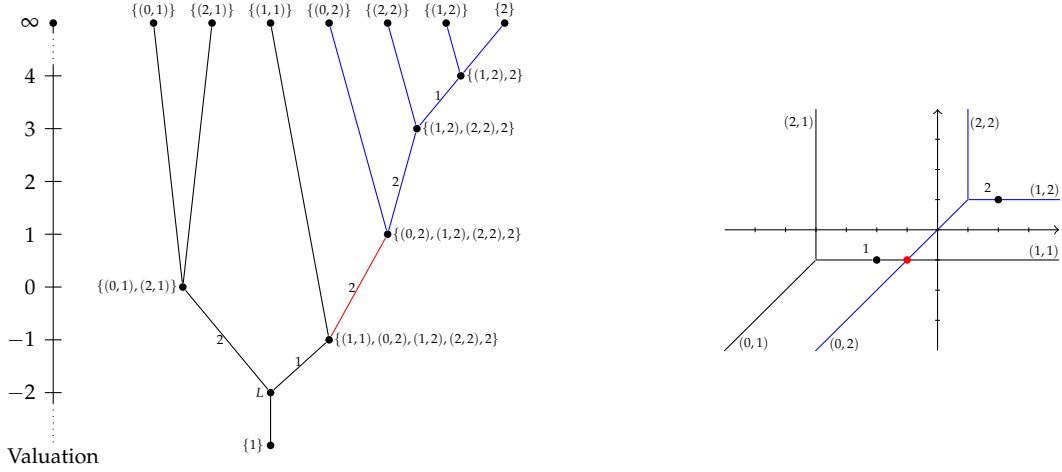
\end{example}

The following proposition will show that the corresponding tropical curve to an algebraic curve (in standard form) corresponds to it via tropicalization.

\begin{proposition}
\label{cProp:Corresponding curve corresponds to tropicalization}
Let $\mathcal C=(\PS 1, (p_\lambda)_{\lambda\in L_0}, \apara)$ be a curve in $\Mlno I$ in standard form. Then the point corresponding to $\mathfrak C_{\mathcal C}$ in $\tMlI I$ is the tropicalization of the point in $\Mlno I$ corresponding to $\mathcal C$.
\end{proposition}

\begin{proof}
We will use the same notation as in Construction \ref{cConstr:Corresponding tropical curve}.  First, we compute the coordinates of $\mathcal C$. By the procedure described at the beginning of this subsection, we need to consider the Plücker coordinates corresponding to the $2\times L_0$ matrix whose $\lambda$-th column is equal to $(0, 1)^t$ if $\lambda=i_0$, and $(1, a_\lambda)^t$ else.  We get Plücker coordinates by taking the $2\times 2$ minors of this matrix. The minor $\amin_{(\lambda,\mu)}$ corresponding to a pair of distinct indices $\lambda,\mu\in L_0$ is equal to $\pm 1$ if $\lambda$ or $\mu$ is equal to $i_0$, and equal to $a_\mu-a_\lambda$ else. The point in $\mathds G_m^{\dPair{L_0}}/\mathds G_m^{L_0}\times \aT r$ corresponding to $\mathcal C$ in the torus embedding of $\Mlno I$ is represented by $\left((\amin_{(\lambda,\mu)})_{(\lambda,\mu)}, (c_i)_{0\leq i\leq r}\right)$. This tropicalizes to the vector represented by $\left(\nu(\amin_{(\lambda,\mu)})_{(\lambda,\mu)},(\nu(c_i))_{0\leq i\leq r}\right)$ in $\mathds R^{\dPair{L_0}}/\Image\Phi\times \tT r$, where $\Phi$ is the morphism $\Phi:\mathds R^ {L_0}\rightarrow \mathds R^{\dPair{L_0}}:(x_\lambda)_\lambda\mapsto (x_\lambda+x_\mu)_{(\lambda,\mu)}$ (see Subsection \ref{cSubSec:Tropical Curves}).

To compute the point in $\mathds R^{\dPair{L_0}}/\Image\Phi\times \tT r$ corresponding to $\mathfrak{C}_{\mathcal C}=(\mathfrak C,(e_\lambda),\tpara)$ we need to compute the distances of the legs. To shorten the notation, we write $v_\lambda\coloneqq v_{e_\lambda}$ for the inner vertex incident to $e_\lambda$. For every pair of distinct indices $\lambda, \mu$ the distance $\dist_{\mathfrak C}(\lambda,\mu)$ is, by definition, equal to the distance between $v_\lambda$ and $v_\mu$. If $\lambda=i_0$, we have $v_\lambda=L$ and thus $\dist_{\mathfrak C}(\lambda,\mu)=\nu(v_\mu)-\nu(L)$. Similarly, we get that $\dist_{\mathfrak C}(\lambda,\mu)=\nu(v_\lambda)-\nu(L)$ if $\mu=i_0$.  Else, let $v_{\{\lambda,\mu\}}$ be the inclusion-minimal vertex containing $\{\lambda,\mu\}$. Then
\begin{equation*}
\dist_{\mathfrak C}(\lambda,\mu)=\nu(v_\lambda)+\nu(v_\mu)-2\nu(v_{\{\lambda,\mu\}})=\nu(v_\lambda)+\nu(v_\mu)-2\nu(a_\mu-a_\lambda)
\end{equation*}
So if we define
\begin{align*}
\tmin_{(\lambda,\mu)}= 
\begin{cases}
\left(\nu(L)-\nu(v_\mu)\right)/2 & \text{, if }\lambda=i_0\\
\left(\nu(L)-\nu(v_\lambda)\right)/2 &\text{, if }\mu=i_0\\
\nu(a_\mu-a_\lambda)-\left(\nu(v_\lambda)+\nu(v_\mu)\right)/2 & \text{, else },
\end{cases}
\end{align*}
the point $\tpl\times\ev {i_0}(\mathfrak C_{\mathcal C})\in\mathds R^{\dPair{L_0}}/\Image\Phi\times \tT r$ corresponding to $\mathfrak C_{\mathcal C}$ is represented by the pair $((\tmin_{(\lambda,\mu)})_{(\lambda,\mu)},(\nu(c_i))_{0\leq i\leq r})$. We may alter the vector $(\tmin_{(\lambda,\mu)})$ by an element in the image of $\Phi$, without changing the represented point. In particular, we can add $\Phi(x_\lambda)$, where we define 
\begin{align*}
x_\lambda= 
\begin{cases}
-\nu(L)/2 & \text{, if }\lambda=i_0\\
\nu(v_\lambda)/2 & \text{, else }.
\end{cases}
\end{align*}
The $(\lambda,\mu)$-th coordinate of the resulting representative of $\tpl(\mathfrak C)$ is equal to $0=\nu(\amin_{(\lambda,\mu)})$, if $\lambda=i_0$ or $\mu=i_0$, and equal to $\nu(a_\mu-a_\lambda)=\nu(\amin_{(\lambda,\mu)})$, else, yielding the desired result.
\end{proof}

Next we will prove that evaluation commutes with tropicalization. This will be the key step in making parametrized rational curves accessible to the methods developed in the next Section.

\begin{proposition}
\label{cProp:Image of corresponding tropical curve is tropicalization of image of curve}
Let $\mathcal C=(\PS 1, (p_\lambda)_{\lambda\in L_0}, \apara)$ be a curve in $\Mlno I$ in standard form with corresponding tropical curve $\mathfrak C_{\mathcal C}=(\mathfrak C,(e_\lambda),\tpara)$. Then $\tpara(\mathfrak C)\cap\VG(\tT r)=\trop(\apara(C)\cap\aT r)$ and for every $i\in I$ we have $\tev i (\mathfrak C_{\mathcal C})=\trop(\ev i (\mathcal C))$.
\end{proposition}

\begin{proof}
Both parts of the statement are instances of the same question, namely how to give a nice description of the tropicalization of a point $\apara(1:a)$ for $a\in K$. We will acquire such a description while proving that $\trop(\apara(C)\cap\aT r)$ is a subset of $\tpara(\mathfrak C)\cap\VG(\tT r)$. So let $x=\trop(\apara(p))$ for some $p\in \PS 1$ not equal to any of the $p_{ij}$. If $p=p_{i_0}=(0:1)$, then $\trop(\apara(p))=\tpara(L)$. Thus we can assume $p=(1:a)$ for some $a\in K$. The image $f(1:a)$ of $p$ in $\PS r$ is represented by
\begin{equation*}
\left(c_0 \prod_{i=1}^d(a-a_{0i}),\dots, c_r\prod_{i=1}^d (a-a_{ri})\right),
\end{equation*}
where we use the same notation as in Construction \ref{cConstr:Corresponding tropical curve}. Therefore, the tropicalization $\trop(\apara(p))$ is represented by
\begin{equation*}
\left(\nu(c_0)+ \sum_{i=1}^d\nu(a-a_{0i}),\dots, \nu(c_r)+\sum_{i=1}^d \nu(a-a_{ri})\right).
\end{equation*}

 Let $r_1<\dots< r_k$ be the elements of $\{\nu(a-a_\lambda)\mid \lambda \in J\}$, and for each $1\leq i\leq k$ define $D_i=\{\lambda\in J \mid \nu(a-a_\lambda)\geq r_i\}$. Then the expression above can be written as
\begin{equation}
\label{cEqu:Tropicalization of a point on the curve}
\sum_{i=0}^r\nu(c_i)s_i+\sum_{i=2}^k (r_i-r_{i-1})s_{D_i}, +\underbrace{r_1 s_{D_1}}_{=0},
\end{equation} 
where $s_{D_1}=0$ because it is represented by $d$-times the all one vector. Let $v_i\in V$ be the inclusion-minimal vertex containing $D_i$. In order to prove that $\trop(\apara(p))$ is the image of a point of $\mathfrak C$ on the path from $L$ to $v_k$ we will collect some properties of the $v_i$. 
First of all, the choice of $v_i$ ensures that $v_i\supseteq v_j$ for $i\leq j$. Moreover, it follows immediately from the definition of $V$ that $\nu(v_i)=\nu(D_i)$ for every $1\leq i\leq k$.  Using the valuation properties of $\nu$ we also see directly that $\nu(D_i)\geq r_i$ for all $i$. We claim that for $i<k$ we even have equality. Assume the opposite, in which case $\nu(D_i)> r_i$. Let $\lambda\in D_i$ such that $\nu(a-a_\lambda)=r_i$, and let $\mu\in D_k$. Then on the one hand we have 
\begin{equation*}
\nu(a_\lambda-a_\mu)\geq \nu(D_i)> r_i=\nu(a-a_\lambda),
\end{equation*}
and hence $\nu(a-a_\mu)=r_i$ by the valuation properties, whereas on the other hand we have $\nu(a-a_\mu)=r_k$ by choice of $\mu$, which obviously contradicts the fact that $r_i\neq r_k$. 

Next we show that $v_i\cap J=D_i$. One inclusion is clear, so let $\lambda\in v_i\cap J$, and choose an arbitrary $\mu\in D_i$. We obtain
\begin{equation*}
\nu(a-a_\lambda)=\nu((a-a_\mu)+(a_\mu-a_\lambda))\geq \min\{\nu(a-a_\mu),\nu(a_\mu-a_\lambda)\}.
\end{equation*}
As $\nu(a-a_\mu)\geq r_i$ by the definition of $D_i$, and $\nu(a_\mu-a_\lambda)\geq \nu(v_i)=\nu(D_i)\geq r_i$, this minimum is greater or equal to $r_i$, proving that $\lambda\in D_i$. 

Now we show that for any $1\leq i < k$, there does not exist a vertex $v\in V$ such that $v_i\supsetneq v\supseteq v_{i+1}$ and $v\cap J \supsetneq v_{i+1}\cap J=D_{i+1}$. We assume the opposite again. Let $\lambda\in v\cap J\setminus v_{i+1}$, and let $\mu\in D_{i+1}$ arbitrary. By choice of $\lambda$ we have $\lambda\in D_i\setminus D_{i+1}$, and hence $\nu(a_\lambda-a)=r_i$. Furthermore, $\nu(a-a_\mu)\geq r_{i+1}>r_i$. Thus, $\nu(a_\lambda-a_\mu)=r_i$ by the valuation properties. But both $\lambda$ and $\mu$ are in $v$, leading to the contradiction
\begin{equation*}
r_i=\nu(a_\lambda-a_\mu)\geq \nu(v)>\nu(v_i)=r_i.
\end{equation*}

Let $L=w_0\supsetneq w_1\supsetneq \dots \supsetneq w_l= v_k$ be the path from $L$ to $v_k$, and for $1\leq i\leq k$ let $j_i$ denote the index such that $w_{j_i}=v_i$. Furthermore, let $s$ be the maximal index such that $\nu(w_s)<r_k$. Since $\nu(v_k)\geq r_k$, and $\nu(v_{k-1})< r_k$ we have $j_{k-1}\leq s< l$ and 
\begin{equation*}
0<r_k-\nu(w_s)\leq \nu(w_{s+1})-\nu(w_s)=\ell(\{w_s,w_{s+1}\}).
\end{equation*}

Using all of this we can rewrite the second summand in Expression \ref{cEqu:Tropicalization of a point on the curve} as
\begin{align*}
\sum_{i=2}^k (r_i-r_{i-1})s_{D_i}&=
\sum_{i=2}^{k-1} (\nu(v_i)-\nu(v_{i-1}))s_{D_i} +(r_k-\nu(v_{k-1}))s_{D_k}\\
&=\sum_{i=2}^{k-1} \sum_{j=j_{i-1}+1}^{j_i}\ell(\{w_{j-1},w_j\})s_{w_i} + \sum_{j=j_{k-1}+1}^{s}\ell(\{w_{j-1},w_j\})s_{w_j}\\ 
&+(r_k-\nu(w_{s}))s_{w_{s+1}}.
\end{align*}
Adding $\sum_i \nu(c_i)s_i$ to both sides we get that 
\begin{equation*}
\trop(\apara(p))= q_{w_s}+(r_k-\nu(w_{s}))\dir(w_s,\{w_s,w_{s+1}\}),
\end{equation*}
and this clearly is the image under $\tpara$ of a point on the edge from $w_s$ to $w_{s+1}$.

We have just seen how to construct a point  $\mathfrak p\in \mathfrak C^\circ$ for a given point $p\in\PS 1$ such that $\trop(\apara(p))=\tpara(\mathfrak p)$. Let us apply this construction to the case when $p=p_\iota$ for some $\iota\in I\setminus\{i_0\}$ before proving the other inclusion $\tpara(\mathfrak C)\cap \VG(\tT r)\subseteq \trop(\apara(C)\cap \aT r)$. With the notations as above, let $\lambda\in D_k$, and let $v\in V$ be the inclusion-minimal vertex containing $\{\iota,\lambda\}$. It has valuation $\nu(v)=\nu(\{\iota,\lambda\})=r_k$, and since $w_{s+1}\cap v\neq\emptyset$ we either have $w_{s+1}\subseteq v$ or $v\subseteq w_{s+1}$. Because $\nu(w_{s+1})\geq r_k$ by construction, we must have $w_{s+1}\subseteq v$. If $v$ strictly contained $w_{s+1}$, then it would also contain $w_s$ and therefore satisfy $r_k=\nu(v)\leq \nu(w_s)<r_k$, a contradiction! We deduce that $v=w_{s+1}$ and hence that $\nu(w_{s+1})=r_k$ and $\iota\in w_{s+1}$. Furthermore, we claim that $w_{s+1}$ is inclusion-minimal with the property of containing $\{\iota\}$ and having nonempty intersection with $J$. Assume the opposite, that is that there exists a vertex $\iota\in u\in V$ with $u\subsetneq w_{s+1}$ and $u\cap J\neq\emptyset$. Then for every $\mu\in u\cap J$ we have $\nu(a_\iota-a_\mu)\geq \nu(u)>\nu(w_{s+1})=r_k$, which is a contradiction by the construction of $r_k$! This proves that 
\begin{equation*}
\tev \iota (\mathfrak C_{\mathcal C})=\tpara(w_{s+1})=\trop(\apara(p_\iota))=\trop(\ev \iota (\mathcal C))
\end{equation*}
for $\iota\in I\setminus\{i_0\}$, and since it is trivially true for $\iota=i_0$ we even get the equality $\tev \iota (\mathfrak C)=\trop(\ev \iota (\mathcal C))$ for all $\iota\in I$.

We continue to prove the equality $\tpara(\mathfrak C)\cap \VG(\tT r)=\trop(\apara(C)\cap \aT r)$. Let $\mathfrak p$ be a point in $\tpara(\mathfrak C)\cap \VG(\tT r)$. Then it is of the form $\mathfrak p=q_v+t\dir(v,\{v,w\})$ for an appropriate edge $\{v\supset w\}$ and some $t\in\VG$ with $0<t\leq\nu(w)-\nu(v)$. If $J\subseteq w$, then $\mathfrak p=\trop(\apara(p_{i_0}))$ and we are done, so assume $w\cap J\subsetneq J$. After replacing $w$ by some vertex on the path from $w$ to $L$ we can also assume that $w\cap J\neq\emptyset$.

In the proof of the other inclusion we have seen that we have to construct an $a\in K$ such that $r\coloneqq\max\{\nu(a-a_\lambda)\mid\lambda\in J\}$ is equal to $\nu(v)+t$, and the set $D\coloneqq \{\lambda\in J\mid \nu(a-a_\lambda)\geq r\}$ is contained in $w$. To do this, choose $\mu\in w\cap J$. Because $\nu(a_\mu-a_\lambda)\geq \nu(v)+t$ for all $\lambda\in w\cap J$, and $\kappa$ is infinite, there exists $\overline a\in K$ such that $\nu(\overline a + a_\mu-a_\lambda)=\nu(v)+t$ for all $\lambda\in w\cap J$. We claim that $a\coloneqq \overline a+a_\mu$ has the desired properties. By construction of $a$, we have 
\begin{equation*}
\nu(a-a_\lambda)=\nu(\overline a+a_\mu-a_\lambda)=\nu(v)+t
\end{equation*}
for all $\lambda\in w\cap J$. On the other hand, if $\lambda\in J\setminus w$, then $\nu(a_\mu-a_\lambda)=\nu(w\cup\{\lambda\})$. This is equal to the valuation of a vertex which strictly contains $w$, and therefore contains $v$. In particular, we have $\nu(a_\mu-a_\lambda)\leq\nu(v)<\nu(v)+t$, and thus 
\begin{equation*}
\nu(a-a_\lambda)=\nu(\overline a+a_\mu-a_\lambda)<\nu(v)+t.
\end{equation*}
As a consequence, we have $r=\nu(v)+t$ and $D=w\cap J\subseteq w$. With what we saw in the proof of the other inclusion this yields $\mathfrak p=\trop(\apara(1:a))$.

\end{proof}

\begin{example}
We want to apply the method described at the beginning of the preceding proof to compute a point $\mathfrak p$ in the corresponding tropical curve of the curve of Example \ref{cEx:Corresponding tropical curve} such that $\tpara(\mathfrak p)=\trop(1:2+t+t^2)$. We have 
\begin{equation*}
\{\nu(2+t+t^2-a_\lambda)\mid \lambda\in J\}=\{-2,-1,1,2\}, 
\end{equation*}
that is in the notation of the preceding proof we have $k=4$ and $r_i$, $D_i$ and $v_i$ as follows:
\begin{center}
{\small\renewcommand{\arraystretch}{1.1}
$\begin{array}{c||c|c|c|c}
i & 		\hphantom{-}1 &		 \hphantom{-}2 &		 3 &		 4 \\
\hline
r_i &   -2 &  	 -1&		1&		2\\
D_i & \hphantom{-}J & \{(1,1),(0,2),(1,2),(2,2)\} & \{(0,2),(1,2),(2,2)\} & \{(1,2),(2,2)\} \\
v_i & L=J\cup\{2\} & D_2\cup\{2\}& D_3\cup\{2\} & D_4\cup\{2\} \\
\nu(v_i) & -2 & -1 & 1 & 3
\end{array}$}
\renewcommand{\arraystretch}{1.0}
\end{center}
Note that we can easily read off the $v_i$ from Figure \ref{cFig:curve}. The path $L=w_0\supsetneq \dots\supsetneq w_l=v_k$ just consists of the $v_i$ in this example, that is $l=3$ and $w_i=v_{i+1}$. The largest index $s$ such that $\nu(w_s)<r_k$ is equal to $2$. We conclude that the construction yields a point on the edge from $v_3$ to $v_4$. More precisely, because $r_k-\nu(v_3)=1=\frac 12 \ell(\{v_3,v_4\})$ the point is exactly in the middle of that edge and, as can be seen in Figure \ref{cFig:curve}, mapped to $(0:0:0)$.
\end{example}

Because every algebraic curve is equivalent to one in standard form, Propositions \ref{cProp:Corresponding curve corresponds to tropicalization} and \ref{cProp:Image of corresponding tropical curve is tropicalization of image of curve} yield the following corollary.

\begin{corollary}
\label{cCor:Tropicalization of curves commutes with everything}
For every $i\in I$ the diagram
\begin{center}
\begin{tikzpicture}[auto]
\matrix[matrix of math nodes,row sep=0.7cm, column sep=0.9cm]{
|(aM)| \Mlno I & |(aT)| \aT r \\
|(tM)| \tMlI I & |(tT)| \tT r\\
};
\begin{scope}[->,font=\footnotesize]
\draw (aM) --node{$\trop$} (tM);
\draw (aT) --node{$\trop$} (tT);
\draw (aM) --node{$\ev i$} (aT);
\draw (tM) --node{$\tev i$} (tT);
\end{scope}
\end{tikzpicture}
\end{center}
is commutative.
\end{corollary}

We finish this section with another result on the evaluation maps.

\begin{proposition}
\label{cProp:Evaluation maps can be extended}
For every $\iota\in I$ the evaluation map $\ev \iota:\Mlno I\rightarrow \aT r$ can be extended to a morphism of tori $\extev \iota:\extMlno \rightarrow \aT r$.
\end{proposition}

\begin{proof}
The morphism of tori
\begin{align*}
\mathds G_m^{\dPair{L_0}}\times\aT r &\rightarrow \qquad\aT r\\
\left((d_{(\lambda,\mu)}),(c_i)_{0\leq i\leq r}\right)&\mapsto \left(c_i\prod_{j=1}^d\frac{d_{\left((i,j),\iota\right)}}{d_{\left((i,j),i_0\right)}}\right)_{0\leq i\leq r}
\end{align*}
obviously induces a morphism $\extev \iota:\extMlno\rightarrow \aT r$. All there is left to show is that it coincides with $\ev \iota$ on $\Mlno I$. So let $\mathcal C=(\PS 1,(p_\lambda),\apara)$ be a representative in standard form of a point in $\Mlno I$. Using the notation of Construction \ref{cConstr:Corresponding tropical curve}, the proof of Proposition \ref{cProp:Corresponding curve corresponds to tropicalization} shows us that $\mathcal C$ is represented by the point $((d_{(\lambda,\mu)}),(c_i))\in \mathds G_m^{\dPair{L_0}}\times\aT r$, where $d_{(\lambda,\mu)}$ is equal to $-1$ if $\lambda=i_0$, equal to $1$ if $\mu=i_0$, and equal to $a_\mu-a_\lambda$ else. Thus $\extev \iota(\mathcal C)$ is represented by $\left(c_i\prod_{j=1}^d (a_\iota-a_{ij})\right)_{0\leq i\leq r}$. This is exactly $\apara(1:a_\iota)$, which is equal to $\ev \iota (\mathcal C)$ by definition.
\end{proof}

\section{Relating Algebraic and Tropical Enumerative Geometry}
\label{tSec:Relating Algebraic and Tropical Enumerative Geometry}

We have seen in the previous section that the problems of counting algebraic curves and tropical curves are very closely related. In fact, given the formulation of the classical problem as the determination of the number of points in intersections of preimages of generic points via the evaluation maps, we obtain the tropical problem by tropicalizing the moduli space and the evaluation maps and then formulating the same problem. The fact that the tropicalization of the algebraic moduli space is again a moduli space in the tropical world is, of course, more than a mere coincidence, but this knowledge is not needed to pose the tropical enumerative problem. So a natural abstraction of the problem of relating an algebraic and tropical enumerative problem consists of
\begin{itemize}
\item A subscheme $M$ of an algebraic torus $T$, and
\item evaluation maps $\ev i: M\rightarrow T_i$ for $1\leq i\leq n$ into algebraic tori $T_i$, which extend to morphisms $\extev i: T\rightarrow T_i$ of algebraic tori.
\end{itemize}
After possibly passing to subtori of the $T_i$ we can assume that the morphisms $\extev i$ are dominant. Note that morphisms of tori are dominant if and only if they are surjective if and only if they are flat. 

Tropicalizing the given algebraic objects we obtain
\begin{itemize}
\item A tropical subvariety $\mathfrak M=\Trop(M)$ of a tropical torus $\ttorus=\Trop(T)$, and
\item evaluation maps $\tev i:\mathfrak M\rightarrow \ttorus_i$ for $1\leq i\leq n$ into tropical tori $\ttorus_i=\Trop(T_i)$, which extend to morphisms $\Trop(\extev i)=\textev i:\ttorus\rightarrow \ttorus_i$ of tropical tori.
\end{itemize}
The fact that the $\extev i$ are dominant implies that the corresponding tropical morphisms $\textev i$ are surjective.

In this situation, we want to prove that for an appropriate number $n$ of generic points $p_i\in T_i$ the number of points in $M$ mapped to $p_i$ by $\ev i$ for all $i$ is finite and equal to the number of points, counted with the appropriate multiplicities, in $\mathfrak M$ mapped to $\mathfrak p_i\coloneqq\trop(p_i)$ by $\tev i$ for all $i$. As we want to define the occurring multiplicities intersection-theoretically, the natural claim is that 
\begin{equation}
\label{tEqu:Numbers are Equal}
|\ev 1^{-1}\{p_1\}\cap\dots\cap\ev n^{-1}\{p_n\}|=\deg(\tev 1^*[\mathfrak p_1]\cdots \tev n^*[\mathfrak p_n]\cdot \mathfrak M)
\end{equation}
is fulfilled for generic $p_i$. By generic we mean that for every point $\mathfrak p=(\mathfrak p_i)\in\VG(\prod_i\ttorus_i)$ there is a subset of $\trop^{-1}\{\mathfrak p\}$ which is dense in $\prod_i T_i$ such that the equation holds for every point $p=(p_i)$ in this set. 

It seems curios that ``counting'' in algebraic geometry translates to ``counting with multiplicities'' in tropical geometry. The reason for this is that even though the intersection $\ev 1^{-1}\{p_1\}\cap\dots\cap\ev n^{-1}\{p_n\}$ is finite and reduced for generic $p_i$, we cannot assure that all of its points tropicalize to distinct points of $\ttorus$. But if several points tropicalize to the same point $\mathfrak p$, it is only natural to count $\mathfrak p$ with multiplicity greater $1$. In fact, we will show that for generic points $p_i$ the multiplicity of any point of $\tev 1^*[\mathfrak p_1]\cdots \tev n^*[\mathfrak p_n]\cdot \mathfrak M$ is exactly the number of points in $\ev 1^{-1}\{p_1\}\cap\dots\cap\ev n^{-1}\{p_n\}$ having it as tropicalization.

 To make the analogy between algebraic and tropical enumerative geometry even more visible, let us assume for a moment that $M$ is the moduli space of parametrized rational curves. In the modern treatment of the subject, the  study of the number on the left hand side of Equation \ref{tEqu:Numbers are Equal} usually uses a modular compactification $\overline M$ of $M$ to which the evaluation maps can be extended. After using a Bertini argument to show that the set $\ev 1^{-1}\{p_1\}\cap\dots\cap\ev n^{-1}\{p_n\}$ is finite, reduced, and contained in $M$ for a generic choice of the $p_i$, one can use intersection theory on $\overline M$ to show that its cardinality is independent of the choice of points and equal to
\begin{equation*}
\deg(\ev 1^*[p_1]\cdots \ev n^*[p_n]\cdot \overline M).
\end{equation*}
This is clearly analogous to the expression on the right hand side of Equation \ref{tEqu:Numbers are Equal}. The only difference between the two expressions, except the use of different fonts, reveals one of the big advantages of passing to the tropical world: in tropical geometry we do not need to compactify the moduli space because intersection products and pull-backs are well-defined as tropical cycles. And even if we allow rational equivalence we can still take degrees of cycle classes. This allows us to move the points into special positions without having to deal with the often very difficult task of finding a suitable modular compactification of the moduli space.

To prove the equality in Equation \ref{tEqu:Numbers are Equal} we will need several auxiliary statements. We start by showing that in the simplest case where $M=T$, and we only have one evaluation map, taking inverse images commutes with tropicalization.

\begin{proposition}
\label{tProp:Pullback and Tropicalization commutes}
Let $f:T\rightarrow T'$ be a dominant morphism of algebraic tori and $\mathfrak f:\ttorus\rightarrow \ttorus'$ its tropicalization. Furthermore, let $X\subseteq T'$ be a pure-dimensional subscheme of $T'$. Then 
\begin{equation*}
\Trop(f^{-1}X)=\mathfrak f^*\Trop(X),
\end{equation*}
that is tropicalization commutes with flat pull-back.
\end{proposition}

\begin{proof}
Let $M$ and $M'$ denote the character lattices of $T$ and $T'$, respectively. Since $f$ is dominant, the induced morphism $f^*:M'\rightarrow M$ is injective and we can assume that $M'$ is a subgroup of $M$. Let $\wM$ denote the saturation of $M'$ in $M$. Then the inclusion $M'\hookrightarrow M$ factors through $\wM$, that is $f^*$ is the composite of the inclusions $M'\hookrightarrow \wM$ and $\wM\hookrightarrow M$. Denoting the torus associated to $\wM$  by $\wT=\Spec(K[\wM])$, we see that we can write $f$ as the composite of two dominant morphism $T\rightarrow\wT$ and $\wT\rightarrow T'$ of which the first one splits and the second one is finite. As both the algebraic and the tropical pull-backs are functorial, it suffices to assume that $f$ is either finite or splits. 

First suppose that $f$ splits. In this case, we can assume that $T=T'\times\bT$ for some torus $\bT$, and $f$ is the projection onto the first coordinate. Then $f^{-1}X=X\times\bT$ and, writing $\Trop(\bT)=\overline{\ttorus}$, it is immediate that the underlying set of its tropicalizations is equal to $|\Trop(X)|\times \overline{\ttorus}$, which is the underlying set of the pull-back of $\Trop(X)$. As for multiplicities, let $\mathfrak t=(\mathfrak t',\overline{\mathfrak t})\in \ttorus\times\overline{\ttorus}$. After an appropriate extension of scalars we may assume that there is a $t=(t',\overline t)\in T'\times\overline T$ that tropicalizes to $\mathfrak t$. The pull-back of $t'^{-1}X$ along $f$ is equal to $t^{-1}(f^{-1}X)$. Together with the flatness of $f$ this implies $\ini_{\mathfrak t}(f^{-1}X)\cong f_\kappa^{-1}\left(\ini_{\mathfrak t'}(X)\right)$ by \cite[Cor. 4.5]{Gub12}, where $f_\kappa$ is the induced morphism $T_\kappa\rightarrow T'_\kappa$ between the initial degenerations of the tori. Because flat pull-backs are compatible with pulling back subschemes (\cite[Lemma 1.7.1]{F98}), the sum of the multiplicities of the components of $\ini_{\mathfrak t}(f^{-1}X)$ is equal to the sum of the multiplicities of the components of $\ini_{\mathfrak t'}(X)$. Since these sums are, by definition, the multiplicities of the tropicalizations of $f^{-1}X$ and $X$ at $\mathfrak t$ and $\mathfrak t'$, respectively, we get the desired result that $\mathfrak f^*\Trop(X)=\Trop(X)\times \overline{\ttorus}=\Trop(f^{-1}X)$ as tropical cycles.

Now suppose that $f$ is finite. Then its degree is equal to the index $[M:M']$. By \cite[Ex. 1.7.4]{F98} the push-forward of the cycle associated to $f^{-1} X$ is equal to $f_*[f^{-1} X]=f_*(f^*[X])=[M:M']\cdot[X]$. Now we tropicalize this equation. As the tropicalization of a torus subscheme is equal to the tropicalization of its associated cycle \cite[Remark 13.12]{Gub12} and tropicalization commutes with push-forward \cite[Thm. 13.17]{Gub12} we get $\mathfrak f_*\Trop(f^{-1} X)=[M:M']\cdot\Trop(X)$. On the other hand, using tropical intersection theory we obtain the equation
\begin{align*}
\mathfrak f_*\mathfrak f^*\Trop(X)=\mathfrak f_*(\mathfrak f^*\Trop(X)\cdot \ttorus)
&= \Trop(X)\cdot  \mathfrak f_*\ttorus=\\&=\Trop(X)\cdot ([M:M']\cdot\ttorus')=[M:M']\cdot\Trop(X),
\end{align*}
where the second equality uses the projection formula \cite[Thm. 8.3]{FR10}. We see that $\Trop(f^{-1}X)$ and $\mathfrak f^*\left(\Trop(X)\right)$ have the same push-forward under $\mathfrak f$. But since $\mathfrak f$ is an isomorphism of $\mathds R$-vector spaces, the push-forward $\mathfrak f_*:Z_*(\ttorus)\rightarrow Z_*(\ttorus')$ is injective. We conclude that $\Trop(f^{-1}X)=\mathfrak f^*\left(\Trop(X)\right)$, finishing the proof.
\end{proof}

The preceding proposition shows that for every point $p\in T_i$ we have the equality $\Trop(\extev i^{-1}\{p\})=\textev i^*[\trop(p)]$. This is helpful, yet we are actually interested in $\ev i$ rather than in $\extev i$. We can recover $\ev i^{-1}\{p\}$ from $\extev i^{-1}\{p\}$ by intersecting with $M$. In a similar vein, $\tev i^*[\trop(p)]\cdot\mathfrak M$ is equal to $\textev i^*[\trop(p)]\cdot \mathfrak M$. So the equality we wish to have is $\Trop(\extev i^{-1}\{p\}\cap M)=\textev i^*[\trop(p)]\cdot\mathfrak M$. Unfortunately, this cannot be expected to be true in general, because intersection and tropicalization do not always commute. However, as shown in \cite[Thm. 5.3.3]{OP13}, we actually get equality if we translate $\extev i^{-1}\{p\}$ by a general $R$-integral point $t\in T$. The translate $t\cdot\extev i^{-1}\{p\}$ is equal to $\extev i^{-1}\{\extev i(t)p\}$, which is the pull-back of a translate of $p$. For this to be of use for us we need that $\extev i(t)$ is a general $R$-integral point of $T_i$, the problem of course being the generality. A positive answer to this is provided by the following proposition.

\begin{proposition}
\label{tProp:Induced morphism of R-integral points is open}
Let $f:T\rightarrow T'$ be a dominant morphism of tori. Then the induced map $f_R:\mathcal T(R) \rightarrow \mathcal T'(R)$ between their respective $R$-integral points is open.
\end{proposition}

\begin{proof}
Let $M$ and $M'$ denote the character lattices of $T$ and $T'$, respectively. Then $\TR$ is canonically in bijection with $\Hom(M, R^*)$ and, analogously, $\TpR$ with $\Hom(M', R^*)$. Identifying $T_\kappa(\kappa)$ and $T'_\kappa(\kappa)$ with $\Hom(M,\kappa^*)$ and $\Hom(M',\kappa^*)$, respectively, we get a commutative diagram
\begin{center}
\begin{tikzpicture}[auto]
\matrix[matrix of math nodes,row sep=0.8cm]{
|(HMR)|\Hom(M,R^*) &[0.9cm] |(HM'R)| \Hom(M',R^*) \\
|(HMk)|\Hom(M,\kappa^*) & |(HM'k)|\Hom(M',\kappa^*),\\
};
\begin{scope}[->,font=\footnotesize]
\draw (HMR)  --node{$f_R$} (HM'R);
\draw (HM'R) --node{$r'$} (HM'k);
\draw (HMR)  --node{$r$} (HMk);
\draw (HMk)  --node{$f_\kappa$} (HM'k);
\end{scope}
\end{tikzpicture}
\end{center}
with the vertical morphisms $r$ and $r'$ being the reduction maps. We recall that the topology on $R$-integral points is defined as the initial topology with respect to the reduction map. Hence, to show the openness it is sufficient to prove that $f_\kappa$ is open and $f_R(r^{-1}U)=(r')^{-1}f_\kappa(U)$ for all open sets $U\subseteq \Hom(M,\kappa^*)$. The first of these assertions follows immediately  from $f_\kappa$ being a dominant morphism of tori over an algebraically closed field. For the latter one it suffices to show that the map from $\Hom(M,R^*)$ into the set-theoretic fiber product of $\Hom(M,\kappa^*)$ and $\Hom(M',R^*)$ over $\Hom(M',\kappa^*)$ is surjective. To show this consider the commutative diagram
\begin{center}
\begin{tikzpicture}[auto]
\matrix[matrix of math nodes,row sep=0.8cm, column sep=0.6cm]{
 |(HM1)|\Hom(M,1+\mathfrak{m})			&	|(HMR)|\Hom(M,R^*) 		&	|(HMk)|\Hom(M,\kappa^*) &	|(2se)|0\\
 |(HM'1)|\Hom(M',1+\mathfrak{m})		&	|(HM'R)| \Hom(M',R^*)	&	|(HM'k)|\Hom(M',\kappa^*)	& |(3se)|0\\
 |(1ze)| 0								&	|(2ze)| 0				&	|(3ze)| 0	.			&\\
};
\begin{scope}[->,font=\footnotesize]
\draw (HM1) -- (HM'1);
\draw (HM'1) -- (1ze);

\draw (HMR)  -- (HM'R);
\draw (HM'R) -- (2ze);

\draw (HMk)  -- (HM'k);
\draw (HM'k) -- (3ze);

\draw (HM1) -- (HMR);
\draw (HMR)  -- (HMk);
\draw (HMk) -- (2se);

\draw (HM'1) -- (HM'R);
\draw (HM'R) -- (HM'k);
\draw (HM'k) -- (3se);
\end{scope}
\end{tikzpicture}
\end{center}
It has exact rows because $M$ and $M'$ are free, and the last two columns are exact since $R^*$ and $\kappa^*$ are divisible. If we can show that the first column is exact as well the rest of the proof is just an easy diagram chase. Thus it suffices to show that $1+\mathfrak m$ is divisible. Let $a\in 1+\mathfrak{m}$ and $n\in \mathds N$. Then the polynomial $x^n-a\in R[x]$ factors in $R[x]$ as $(x-b_1)\cdots (x-b_n)$ for appropriate $b_i\in R$, because $K$ is algebraically closed and $R$ is integrally closed	 in $K$.  Reducing modulo $\mathfrak{m}$ yields $x^n-1=(x-b_1)\cdots (x-b_n)$ in $\kappa[x]$. Plugging in  $1$ we get that $b_i\in 1+\mathfrak{m}$ for some $i$. Hence $a$ has an $n$-th root in $1+\mathfrak{m}$, proving the desired divisibility.
\end{proof}

We have seen so far that if $p$ is a point in $T_i$ and $t\in T_i$ is a general element with tropicalization $0$, then $\Trop(\ev i^{-1}\{t\cdot p\})=\tev i^*[\trop(p)]\cdot\mathfrak M$, that is pull-backs of general points under an evaluation map commute with tropicalization. Looking back at Equation \ref{tEqu:Numbers are Equal} this is exactly what we need, only that we need to consider several evaluation maps at once. But with the result of the following lemma, the general situation can be handled analogously.

\begin{lemma}
\label{tLem:Generic Intersection and Tropicalization commute}
Let $X_1,\dots, X_n$ be subschemes of the torus $T$. Then for general $R$-integral points $(t_2,\dots, t_n)$ of $T^{n-1}$ we have 
\begin{equation*}
\Trop\left(X_1\cap t_2X_2\cap\dots\cap t_nX_n\right)=\Trop(X_1)\cdots\Trop(X_n).
\end{equation*}
\end{lemma}

\begin{proof}
For $n=2$ this is just the statement of \cite[Thm. 5.3.3]{OP13}. To prove the general case let $\Delta:T\rightarrow T^n$ be the diagonal embedding, which makes $T$ a subscheme of $T^n$. It corresponds to the diagonal embedding $\mathfrak D:\ttorus\rightarrow \ttorus^n$ on the tropical side. By the $n=2$ case, for every general $R$-integral point $t=(t_1,\dots, t_n)\in T^n$ we have 
\begin{equation*}
\Trop(\Delta\cap t(X_1\times\dots\times X_n))=\Trop(\Delta)\cdot\Trop\left(X_1\times\dots\times X_n\right). 
\end{equation*}
Since $\Delta\cap t(X_1\times\dots\times X_n)$ is the push-forward of $t_1 X_1\cap\dots\cap t_n X_n$, and tropicalization commutes with push-forwards, the term on the left hand side is equal to $\mathfrak D_*(\Trop(t_1 X_1\cap \dots\cap t_n X_n))$. On the other hand, by  \cite[Prop. 5.2.2]{OP13} the term on the right is equal to $\mathfrak D_*(\Trop(X_1)\cdots\Trop(X_n))$. Because $\mathfrak D_*$ is one-to-one, it follows that $\Trop(t_1 X_1\cap \dots\cap t_n X_n)=\Trop(X_1)\cdots\Trop(X_n)$. After multiplying by $t_1^{-1}$ we can assume that $t_1=1$, and Proposition \ref{tProp:Induced morphism of R-integral points is open} ensures that the resulting $R$-rational point $(t_2,\dots,t_n)\in T^{n-1}$ still is generic.
\end{proof}

\begin{theorem}
\label{tThm:Tropicalization commutes with pullback of generic translate}
Using the notation introduced at the beginning of this section, let $P_i\subseteq T_i$ be a pure-dimensional closed subscheme of $T_i$ with tropicalization $\Trop(P_i)=\mathfrak P_i$. Then for generic $R$-integral points $(t_1,\dots, t_k)\in\prod_i T_i$ we have 
\begin{equation*}
\label{tEqu:Tropicalization and generic preimage commute}
\Trop(\ev 1^{-1}(t_1P_1)\cap\dots\cap\ev n^{-1}(t_n P_n))=\tev 1^*(\mathfrak P_1)\cdots \tev n^*(\mathfrak P_n)\cdot\mathfrak M.
\end{equation*}
In particular, if $\mathfrak p=(\mathfrak p_1,\dots,\mathfrak p_n)\in\VG(\prod_i\ttorus_i)$, then for generic  points $p=(p_1,\dots,p_n)$ of $\trop^{-1}\{\mathfrak p\}$ we have
\begin{equation*}
\Trop\left(\ev 1^{-1}\{p_1\}\cap\dots\cap\ev n^{-1}\{p_n\}\right)=\tev 1^*[\mathfrak p_1]\cdots \tev n^*[\mathfrak p_n]\cdot\mathfrak M.
\end{equation*}
\end{theorem}

\begin{proof}
Proposition \ref{tProp:Pullback and Tropicalization commutes} and Lemma \ref{tLem:Generic Intersection and Tropicalization commute} tell us that for a generic $R$-integral point $(s_1,\dots,s_n)\in T^n$ we have 
\begin{equation}
\label{tEqu:Tropicalization and generic preimage commute. Aux1}
\Trop\left(s_1 \cdot\extev 1^{-1}(P_1)\cap\dots\cap s_n\cdot\extev n^{-1}(P_n)\cap M\right)
=\textev 1^*(\mathfrak P_1)\cdots\textev n^*(\mathfrak P_n)\cdot\mathfrak M.
\end{equation}
Let $t=(t_1,\dots,t_n)=(\extev 1(s_1),\dots,\extev n(s_n))\in\prod_i T_i$. Because all evaluation maps are dominant we can apply Proposition \ref{tProp:Induced morphism of R-integral points is open} and see that $t$ is again generic. Furthermore, the left-hand side of Equation \ref{tEqu:Tropicalization and generic preimage commute. Aux1} is equal to 
\begin{equation*}
\Trop\left(\extev 1^{-1}(t_1P_1)\cap\dots\cap \extev n^{-1}(t_nP_n)\cap M\right) 
=\Trop\left(\ev 1^{-1}(t_1P_1)\cap\dots\cap\ev n^{-1}(t_n P_n)\right).
\end{equation*}
As the right-hand side of Equation \ref{tEqu:Tropicalization and generic preimage commute. Aux1} is equal to $\tev 1^*(\mathfrak P_1)\cdots\tev n^*(\mathfrak P_n)\cdot\mathfrak M$ this proves the main statement of the theorem. 

For the ``in particular'' statement let $\mathfrak p=(\mathfrak p_1,\dots,\mathfrak p_n)\in\VG(\prod_i\ttorus_i)$, and choose an arbitrary $s=(s_1,\dots,s_n)\in\trop^{-1}\{\mathfrak p\}$. By what we already proved, we have 
\begin{equation*}
\Trop\left(\ev 1^{-1}\{t_1 s_1\}\cap\dots\cap\ev n^{-1}\{t_n s_n\}\right)=\tev 1^*[\mathfrak p_1]\cdots \tev n^*[\mathfrak p_n]\cdot\mathfrak M.
\end{equation*}
for generic $t=(t_1,\dots, t_n)\in\trop^{-1}\{0\}$. Since $\trop^{-1}\{\mathfrak p\}$ is a $\trop^{-1}\{0\}$-torsor, the elements of the form $ts$ for generic $R$-rational points $t$ are generic in $\trop^{-1}\{\mathfrak p\}$.
\end{proof}

As a corollary to Theorem \ref{tThm:Tropicalization commutes with pullback of generic translate} we obtain the result that Equation \ref{tEqu:Numbers are Equal} is indeed generally true.

\begin{corollary}
\label{tCor:Tropical and Algebraic Numbers are Equal}
With the notation as in the beginning of this section, assume that $M$ is integral, and $\dim(M)=\sum_i \dim(T_i)$. Then for every $\mathfrak p=(\mathfrak p_i)\in\VG(\prod_i \ttorus_i)$ there exists $p=(p_i)\in\prod_i T_i$ with $\trop(p)=\mathfrak p$ such that 
\begin{equation*}
|\ev 1^{-1}\{p_1\}\cap\dots\cap\ev n^{-1}\{p_n\}|=\deg(\tev 1^*[\mathfrak p_1]\cdots \tev n^*[\mathfrak p_n]\cdot \mathfrak M),
\end{equation*}
and such that the multiplicity of every point $\mathfrak q\in \left|\tev 1^*[\mathfrak p_1]\cdots \tev n^*[\mathfrak p_n]\cdot \mathfrak M\right|$  is equal to the number of points in $\ev 1^{-1}\{p_1\}\cap\dots\cap\ev n^{-1}\{p_n\}$ tropicalizing to it. Furthermore, the set of such $p$ is dense in $\prod_i T_i$. More precisely, $p$ can be any point in the intersection of a nonempty open subset of $\trop^{-1}\{\mathfrak p\}$ and a nonempty open subset of $\prod_i T_i$.
\end{corollary}

\begin{proof}
Let $\mathfrak p=(\mathfrak p_i)\in\prod_i \ttorus_i$. Then by Theorem \ref{tThm:Tropicalization commutes with pullback of generic translate} there is a nonempty open subset $W\subseteq \trop^{-1}\{\mathfrak q\}$ such that 
\begin{equation*}
\Trop(\ev 1^{-1}\{p_1\}\cap\dots\cap\ev n^{-1}\{p_n\})=\tev 1^*[\mathfrak p_1]\cdots \tev n^*[\mathfrak p_n]\cdot\mathfrak M
\end{equation*}
for every $(p_1,\dots,p_n)\in W$. Using a Bertini-type argument, we find a dense open subset $U\subseteq\prod_i T_i$ such that $\ev 1^{-1}\{p_1\}\cap\dots\cap\ev n^{-1}\{p_n\}$ is a reduced finite subscheme of $M$. By \cite[Thm. 4.2.5]{OP13}, every nonempty open subset of a fiber of the tropicalization map $\trop:\prod T_i\rightarrow \prod\ttorus_i $ is dense in $\prod T_i$. Hence $W\cap U$ is dense in $\prod T_i$, as well. The corollary follows, since the tropicalization of a reduced finite subscheme of a torus is nothing but the union of the tropicalizations of its points where each point has as multiplicity the number of points tropicalizing to it.
\end{proof}

\section{Applications}
\label{aSec:Applications}

In this section we want to give applications of our results of the previous section. We will start by revisiting our motivational example of labeled parametrized curves in projective space. Because we already studied this case in detail in Section \ref{cSec:Curves}, and our correspondence theorem was modeled after them, we immediately obtain the equality of their algebraic and tropical enumerative numbers. However, it is rather restrictive only to consider curves in projective space on the algebraic side and degree $d$ curves in $\tT r$ on the tropical side. Therefore, we generalize our results of Section \ref{cSec:Curves} to the case of curves in arbitrary normal toric varieties without torus factors in algebraic geometry, and curves of arbitrary degree in arbitrary tropical tori in tropical geometry. Finally, we will give an application to the problem of counting tropical lines in smooth tropical cubic surfaces and show that it is possible to give a purely tropical definition of multiplicities on these lines which reflects their relative realizability. In particular, the multiplicities of the lines on a given smooth cubic surface add up to $27$.

\subsection{Rational Curves in $\PS r$ revisited}

Let $K$ be an algebraically closed valued field with characteristic $0$, and let $d$, $r$ and $n$ be positive integers such that $nr= n+r+rd+d-3$. Take $M=\Mln$, and for $\ev i $ the usual evaluation maps. Each of the evaluation maps is dominant as we can freely choose where one of the marked points is mapped to. So by the choice of $n$, $r$, and $d$, all requirements of Corollary \ref{tCor:Tropical and Algebraic Numbers are Equal} are fulfilled. We see that if $\mathfrak p=(\mathfrak p_i)_i\in\VG(\prod_i\tT r)$ and $p=(p_i)_i\in \prod_i\aT r$ is generic and tropicalizes to $\mathfrak p$, then 
\begin{equation}
\label{aEqu:Numbers are equal for rational curves in projective space}
|\ev 1^{-1}\{p_1\}\cap\dots\cap\ev n^{-1}\{p_n\}|=\deg(\tev 1^*[\mathfrak p_1]\cdots \tev n^*[\mathfrak p_n]\cdot \mathfrak M),
\end{equation}
where $\mathfrak M=\Trop(M)$ and $\tev i= \Trop(\ev i)$. By Corollary \ref{cCor:Tropicalization of curves commutes with everything} we have $\mathfrak M=\tMln$ and the $\tev i$ are just the tropical evaluation maps. Thus Equation \ref{aEqu:Numbers are equal for rational curves in projective space} states that the number of labeled parametrized rational algebraic curves of degree $d$ through the points $p_i$ is equal to the number of labeled parametrized rational tropical curves of degree $d$ through the points $\mathfrak p_i$, counted with multiplicities. Moreover, the multiplicity of a tropical curve $\mathfrak C$ through the $\mathfrak p_i$ is equal to the number of algebraic curves through the $p_i$ which tropicalize to $\mathfrak C$. In the algebraic as well as in the tropical situation the number obtained in this way is not equal to the number of (unlabeled) rational curves through the given points, because we can label the same curve in different ways. But in both cases the number we are interested in can be obtained from the result in the same way, namely by dividing by $(d!)^{(r+1)}$. So to sum up, the number of rational curves in $\PS r$ of certain degree through an appropriate number of points is equal to the number of rational tropical curves in $\tT r$ of the same degree through the tropicalizations of these points, provided that we count them with the correct multiplicities. This reproves the correspondence theorem of \cite{NS06} in the special case of curves in $\PS r$. 

In tropical geometry it is easy to see that the number of curves through given points is independent from the choice of points. This is because taking pull-backs and taking degrees both respect rational equivalence. So with our correspondence theorem we also see that the number of algebraic curves through given points is independent from the choice of the points, as long as they are generic. We have to be careful though, because in our situation being generic always means that it holds for all points in open subsets of the fibers of the tropicalization map $\trop:\prod_i \aT r\rightarrow \prod_i \tT r$. These sets are all dense in $\prod_i \aT r$, yet this notion of generality differs from the usual one which demands that the points can be chosen in some dense open subset of $\prod_i \aT r$. But if we equip $K$ with the trivial valuation these two notions of generality coincide and we get that for $n$ generic points in $\aT r$ (in the usual sense) the number of rational curves of degree $d$ through them is equal to
\begin{equation*}
(d!)^{-r-1}\deg\left(\tev 1^*[0]\dots\tev n^*[0]\cdot \tMlno nrd\right).
\end{equation*}

We see that if we are given an algebraic enumerative problem as described in Section \ref{tSec:Relating Algebraic and Tropical Enumerative Geometry} over an arbitrary field $K$ (algebraically closed and of characteristic $0$, of course) we can equip $K$ with the trivial valuation and apply  Corollary \ref{tCor:Tropical and Algebraic Numbers are Equal}. In this way we obtain the independence of the enumerative number in question of the choice of generic points, and also a way to compute this number tropically. This relieves us from the usually very difficult task of compactifying the algebraic moduli space in a suitable manner. Knowing that the enumerative information is completely encoded in the tropicalizations, one can then proceed to compute the tropical enumerative numbers and get back the algebraic ones. That reducing an algebraic enumerative problem to a tropical one really makes the problem easier is, of course, not clear a priori. But results like Mikhalkin's lattice path algorithm to compute enumerative invariants of toric surfaces in \cite{Mik05}, or the purely tropical proof of Kontsevich's formula in \cite{GM08}, have shown that the chances are good.

Corollary \ref{tCor:Tropical and Algebraic Numbers are Equal} is only applicable for pull-backs of points, yet sometimes one wishes to pull-back higher dimensional cycles as well. Under certain conditions we can apply Theorem \ref{tThm:Tropicalization commutes with pullback of generic translate} to this situation as well. We illustrate this at the example of counting lines with nonempty intersection with four given lines in $\PS 3$. Let $\overline L_1,\dots, \overline L_4$ be four generic lines in $\PS 3$. As they are generic, they are not contained in any coordinate hyperplane. Define $L_i\coloneqq \overline L_i\cap \aT 3$. Again, as the chosen lines are generic, no line passes three of them and intersects the fourth one in a point contained in one of the coordinate hyperplanes. Therefore, the number of lines passing through all four lines can be directly deduced from the number of elements in the intersection $\ev 1^{-1}(L_1)\cap\dots\cap \ev 4^{-1}(L_4)\subseteq\MlIrd 431$. Applying Theorem \ref{tThm:Tropicalization commutes with pullback of generic translate} we get that for generic $R$-integral points $t_i\in\aT 3$ we have
\begin{equation*}
\Trop\left(\ev 1^{-1}(t_1L_1)\cap\dots\cap\ev 4^{-1}(t_4L_4)\right)
=\tev 1^* (\mathfrak L_1)\cdots\tev 4 ^*(\mathfrak L_4)\cdot \tMlno 431,
\end{equation*}
where $\mathfrak L_i=\Trop(L_i)$. Proceeding as in the proof of Corollary \ref{tCor:Tropical and Algebraic Numbers are Equal}, a Bertini-type argument shows that for generic $t_i$ the subscheme $\ev 1^{-1}(t_1L_1)\cap\dots\cap \ev 4^{-1}(t_4L_4)$ of $\MlIrd 431$ is reduced and finite, hence we get the desired equality
\begin{equation*}
\left|\ev 1^{-1}(t_1L_1)\cap\dots\cap\ev 4^{-1}(t_4L_4)\right|=\deg\left(\tev 1^*(\mathfrak L_1)\cdots\tev 4^* (\mathfrak L_4)\cdot\tMlno 431\right).
\end{equation*}
This number seems to depend on the tropicalization of the four lines. However, since we are only interested in general lines, we can assume that all four of them intersect the four coordinate hyperplanes in four respective points. In this case, they all tropicalize to tropical lines consisting of four rays in standard directions and possibly one bounded edge. All those lines are rationally equivalent, proving the independence of the choice of generic lines. Again, one needs to be careful what generic means, but one can always assume the valuation to be trivial to obtain the usual sense of generality.
Another important issue is that in order to obtain information about the number of lines passing through four general lines, it is essential that $t_i L_i$ is a line again. That is, we need that the type of subscheme we want to pull back is translation invariant.

\subsection{General Tropical Degrees and Relative Invariants}
\label{aSubSec:GeneralDegrees}

So far we only considered curves in projective space intersecting the coordinate hyperplanes transversally, which correspond to tropical curves in $\tT r$ whose unbounded rays all have standard directions. It is natural to ask if we can generalize our results and also allow curves in general toric varieties intersecting the boundary divisors with multiplicities on the algebraic side, and tropical curves with arbitrary degrees on the tropical side. As we will see shortly, our results from Section \ref{cSec:Curves} can be modified in a way that gives an affirmative answer to this question.

Let $\Sigma$ be a (rational, polyhedral, strongly convex) fan in a tropical torus $\ttorus$ which is spanned by the rays of $\Sigma$. On the algebraic side we will consider rational curves in the toric variety $X=X_\Sigma$ associated to $\Sigma$, which intersect all torus orbits in the expected dimension, whereas on the tropical side we will consider curves whose unbounded rays are, after translating their initial points to the origin, contained in $\Sigma(1)$, the set of rays of $\Sigma$. As in the case of projective varieties we consider labeled curves, that is for a curve $\PS 1\rightarrow X$ we will mark the preimages of the boundary divisors $D_\rho$ for rays $\rho\in \Sigma(1)$. But in contrast to our previous approach we will allow our curves to intersect the boundary divisors at the labeled points with fixed multiplicity, possibly not equal to $1$.

Because the curves under consideration do not intersect the torus orbits of codimension greater $1$, we may assume that $\Sigma$ only consists of its rays and $\{0\}$. In particular, $X$ is smooth and each $D_\rho$ is a Cartier divisor. Hence, for every parametrized curve $\mathcal C=(\apara:\PS 1\rightarrow X)$ there are pull-backs $\apara^*(D_\rho)$ for all $\rho\in\Sigma(1)$. Defining $d_\rho\coloneqq\deg(\apara^*(D_\rho))$, the relations of the $D_\rho$ in $\Pic(X)$ entail that $\sum_\rho d_\rho u_\rho=0$, where $u_\rho=u_{\rho/0}$ denotes the primitive lattice vector of $\rho$. We call $d\coloneqq(d_\rho)_\rho$ the degree of $\mathcal C$. If we label the intersections of $\mathcal C$ with the boundary divisors, that is the elements of $\bigcup_\rho \apara^*(D_\rho)$, with elements in a finite set $J$, we obtain as additional data a map $\pi:J\rightarrow \Sigma(1)$ indicating on which divisors the labeled points lie, and a map $\omega:J\rightarrow \mathds N_{>0}$ indicating the multiplicity with which the point corresponding to $j\in J$ intersects $D_{\pi(j)}$. The degree $d$ can be reconstructed from the triple $(J,\pi,\omega)$ via the identity
\begin{equation*}
d_\rho(J,\pi,\omega)\coloneqq \sum_{j\in\pi^{-1}\{\rho\}}\omega(j)=d_\rho.
\end{equation*}
This motivates the definition of a \emph{tangency condition} as a triple $c=(J,\pi,\omega)$, where $J$ is a finite set and $\pi$ and $\omega$ are maps as above, such that 
\begin{equation*}
\sum_\rho d_\rho(c)u_\rho=0.
\end{equation*}
The family $d(c)\coloneqq (d_\rho(c))_\rho$ will be called the \emph{degree associated to c}. Having established this notion we can define the objects we want to consider:

\begin{definition}
Let $c=(J,\pi,\omega)$ be a tangency condition and $I$ a finite set disjoint from $J$. As in Section \ref{cSec:Curves} we set $L_0\coloneqq I\cup J$. A \emph{labeled parametrized $I$-marked (rational) curve in $X$ satisfying $c$} is defined as a triple $(C,(p_\lambda)_{\lambda\in L_0},\apara)$ consisting of an abstract rational curve $C$, a family $(p_\lambda)_\lambda$ of distinct points of $C$, and a morphism $\apara:C\rightarrow X$ whose image is not contained in the boundary of $X$, and such that
\begin{equation*}
\apara^*(D_\rho)=\sum_{\lambda\in\pi^{-1}\{\rho\}} \omega(\lambda)p_\lambda
\end{equation*}
for all $\rho\in\Sigma(1)$.
\end{definition}

We call two labeled parametrized $I$-marked curves isomorphic if there is an isomorphism between their underlying abstract curves respecting the marks, labels, and the morphisms, and denote by $\MlIrd I{X}c$ the set of isomorphism classes of labeled parametrized $I$-marked curves in $X$ satisfying $c$. In the following we will show that $\MlIrd I{X}c$ can be tropicalized in a way analogous to $\Mlno I$. First, let us see how to embed $\MlIrd I{X}c$ into an algebraic torus. In the case of curves in projective space we have used that every map $\PS 1\rightarrow \PS r$ of degree $d$ is given by $(r+1)$ homogeneous polynomials of degree $d$. Now that we work with more general toric varieties we can use the generalization of this result  \cite[Thm. 2.1]{Cox95} which says that every morphism $\PS 1\rightarrow X$ of degree $d(c)$ is given by a family of homogeneous polynomials $(f_\rho)_{\rho\in\Sigma(1)}$ such that $f_\rho$ has degree $d_\rho(c)$ for all $\rho\in\Sigma(1)$. The condition that the $(r+1)$ homogeneous polynomials in the case of morphisms into $\PS r$ must not vanish simultaneously, has its analogy in the condition that for $p\in \PS 1$ the point $(f_\rho(p))_\rho\in\mathds A^{\Sigma(1)}$ must not lie in the subvariety $Z$ of $\mathds A^{\Sigma(1)}$ given by the equations $\prod_{\rho'\nleq\rho}x_{\rho'}$ for $\rho\in\Sigma(1)$ (Remember that the rays of $\Sigma$ are also its maximal cones by assumption). The morphism $\PS 1\rightarrow X$ can then be recovered from the $(f_\rho)$ by considering the natural identification of $X$ with the quotient of $ U\coloneqq \mathds A^{\Sigma(1)}\setminus Z$ by the algebraic torus $G\coloneqq\Hom_{\mathds Z}(\Pic(X),\mathds G_m)$ \cite[Thm. 5.1.11]{CLS}. Let us fix some $i_0\in I$. Since, analogously to the $X=\PS r$ case, two collections of homogeneous polynomials define the same morphism to $X$ if and only if they differ by the action of an element of $G$, we see that a curve $(C,(p_\lambda),\apara)$ is uniquely determined by the points $(p_\lambda)$, and the image $\apara(p_{i_0})\in T$ in the big open torus $T$ of $X$. This gives us an embedding of $\MlIrd I{X}c$ into the algebraic torus $\mathds G_m^{\dPair {L_0}}/\mathds G_m^{L_0}\times T$ similar to the torus embedding of $\Mlno I$ presented in Subsection \ref{cSubSec:Tropicalizing}.

The tangency condition $c$ determines a tropical degree $\Delta=\Delta(c)$ in $\ttorus$ with index set $J$, namely by assigning to $j\in J$ the integral vector $\omega(j)u_{\pi(j)}\in\nt$. The same argument as in Subsection \ref{cSubSec:Tropicalizing} shows that the tropicalization of $\MlIrd I{X}c$ is equal to the tropical moduli space $\tMlno I{\ttorus}{\Delta}$. The reader will not be surprised to learn that all constructions and results from Subsection \ref{cSubSec:Tropicalizing} can easily be carried over to the more general situation. In fact, looking back at the construction of the corresponding tropical curve to a labeled parametrized curve in standard form (Construction \ref{cConstr:Corresponding tropical curve}), we see that all we need to change is the definition of the vectors $s_A$ in Equation \ref{cEqu:Definition of s_a}, namely by defining 
\begin{equation*}
s_A\coloneqq \sum_{j\in A\cap J}\Delta(j).
\end{equation*}
Like this we get a corresponding tropical curve $\mathfrak C_{\mathcal C}$ to every labeled parametrized $I$-marked curve $\mathcal C=(\PS 1,(p_\lambda),\apara)$ satisfying $c$ and being in standard form. The remaining results of Subsection \ref{cSubSec:Tropicalizing} can then be copied almost literally to prove the analogous results in the more general situation. In particular, we get 

\pagebreak
\begin{corollary}
\label{aCor:Tropicalization of curves commutes with everything}
For every $i\in I$ the diagram
\begin{center}
\begin{tikzpicture}[auto]
\matrix[matrix of math nodes,row sep=0.7cm, column sep=0.9cm]{
|(aM)| \MlIrd I{X}c & |(aT)| T \\
|(tM)| \tMlno I{\ttorus}{\Delta} & |(tT)| \ttorus\\
};
\begin{scope}[->,font=\footnotesize]
\draw (aM) --node{$\trop$} (tM);
\draw (aT) --node{$\trop$} (tT);
\draw (aM) --node{$\ev i$} (aT);
\draw (tM) --node{$\tev i$} (tT);
\end{scope}
\end{tikzpicture}
\end{center}
is commutative. Furthermore, the evaluation morphisms $\ev i$ can be extended to morphisms $\extev i:\mathds G^{\dPair {L_0}}_m/\mathds G_m^{L_0}\times T\rightarrow T$ of algebraic tori.
\end{corollary}

So far we have only considered the evaluation maps at the marked points, which are mapped to the big torus $T$ of $X$. However, in some applications, for example when studying relative invariants, one also wishes to evaluate at the labeled points. This can be done algebraically as well as tropically. In algebraic geometry, for every $j\in J$ and curve $(C,(p_\lambda),\apara)\in \MlIrd I{X}c$ the image $\apara(p_j)$ is a well-defined point in the torus orbit $O(j)\coloneqq O(\pi(j))$ corresponding to $\pi(j)\in\Sigma(1)$. This yields a map
\begin{equation*}
\ev j: \MlIrd I{X}c\rightarrow O(j),
\end{equation*}
the evaluation map at $j$. As in the case of evaluations at marked points, these evaluation maps are restrictions of morphisms of tori:

\begin{proposition}
\label{aProp:Evaluation maps can be extended}
For every $j\in J$ the map $\ev j:\MlIrd I{X}c\rightarrow O(j)$ can be extended to a morphism $\extev j:\mathds G^{\dPair {L_0}}_m/\mathds G_m^{L_0}\times T\rightarrow O(j)$ of algebraic tori. In particular, the evaluation maps are morphisms.
\end{proposition}

\begin{proof}
The natural toric morphism $U=\mathds A^{\Sigma(1)}\setminus Z\rightarrow X$ induces a torus morphism from the torus orbit $O_U(j)$ of $U$ corresponding to $\pi(j)$ (or rather to the cone generated by the basis vector of $\mathds A^{\Sigma(1)}$ corresponding to $\pi(j)$)  to the torus orbit $O(j)$ of $X$. The orbit $O_U(j)$ can naturally be identified with $\mathds G_m^{\Sigma(1)\setminus\{\pi(j)\}}$ and therefore it makes sense to define the morphism of algebraic tori
\begin{align*}
\mathds G^{\dPair{L_0}}_m\times \mathds G_m^{\Sigma(1)} &\rightarrow\quad O_U(j)\\
\left((d_{(\lambda,\mu)}),(c_\rho)\right)
&\mapsto \left (c_\rho \prod_{\lambda\in\pi^{-1}\{\rho\}}\left(\frac{d_{(\lambda,j)}}{d_{(\lambda,i_0)}}\right)^{\omega(\lambda)}\right)_{\rho\neq \pi(j)}.
\end{align*}
Taking the composite of this morphism with the morphism $O_U(j)\rightarrow O(j)$ we obtain a morphism $\mathds G^{\dPair{L_0}}_m\times \mathds G_m^{\Sigma(1)}\rightarrow O(j)$ of algebraic tori which obviously contains the naturally embedded $\mathds G_m^{L_0}\times G$ (remember that $G=\Hom_{\mathds Z}(\Pic(X),\mathds G_m)$) in its kernel. Hence it induces a morphism
\begin{equation*}
\extev j: \mathds G^{\dPair{L_0}}_m/\mathds G_m^{L_0}\times T\rightarrow O(j)
\end{equation*}
of algebraic tori. That this morphism coincides with $\ev j$ on $\MlIrd I{X}c$ can be seen in the same way as in Proposition \ref{cProp:Evaluation maps can be extended}.
\end{proof}

In the tropical world the situation is similar. Given $j\in J$ and a tropical curve $(\mathfrak C,(e_\lambda),\tpara)\in \tMlno I{\ttorus}{\Delta}$, the images of the points of $e_j$ all differ by an element in the span of $\pi(j)$. Hence we get a well-defined point in $\tOrb(j)=\tOrb(\pi(j))\coloneqq\ttorus/\Lin(\pi(j))$, yielding a map
\begin{equation*}
\tev j:\tMlno I{\ttorus}{\Delta}\rightarrow \tOrb(j),
\end{equation*}
the evaluation map at $j$. The space $\tOrb(j)$ is a tropical torus, its lattice being the image of $\nt$ under the projection map $\ttorus\rightarrow \tOrb(j)$. Because the lattice $N_{\tOrb(j)}$ can naturally be identified with the lattice of 1-psgs of $O(j)$, the tropical torus $\tOrb(j)$ is the tropicalization of $O(j)$. Therefore, we obtain a diagram analogous to the one in Corollary \ref{aCor:Tropicalization of curves commutes with everything}, which also is commutative as the following proposition shows.

\begin{proposition}
For every $j\in J$ the diagram 
\begin{center}
\begin{tikzpicture}[auto]
\matrix[matrix of math nodes,row sep=0.7cm, column sep=0.9cm]{
|(aM)| \MlIrd I{X}c & |(aT)| O(j) \\
|(tM)| \tMlno I{\ttorus}{\Delta} & |(tT)| \tOrb(j)\\
};
\begin{scope}[->,font=\footnotesize]
\draw (aM) --node{$\trop$} (tM);
\draw (aT) --node{$\trop$} (tT);
\draw (aM) --node{$\ev j$} (aT);
\draw (tM) --node{$\tev j$} (tT);
\end{scope}
\end{tikzpicture}
\end{center}
is commutative.
\end{proposition}

\begin{proof}
Let $\mathcal C=(\PS 1, (p_\lambda), \apara)$ be in standard form, and use the notation of Construction \ref{cConstr:Corresponding tropical curve}. The image of the point $p_j$ under $\apara$ is represented by the point
\begin{equation*}
\left(c_\rho\prod_{\lambda\in\pi^{-1}\{\rho\}}(a_j-a_\lambda)^{\omega(\lambda)}\right)_{\rho\neq\pi(j)}
\end{equation*}
in the torus orbit $\mathds G_m^{\Sigma(1)\setminus\{\pi(j)\}}$ of $U$ corresponding to $\pi(j)$. Therefore, $\trop\left(\ev j(\mathcal C)\right)$ is represented by 
\begin{equation}
\label{aEqu:Representative of Tropicalization}
\sum_{\rho\in\Sigma(1)}\nu(c_\rho)u_\rho +\sum_{\lambda\in J\setminus\{j\}}\nu(a_j-a_\lambda)\Delta(\lambda).
\end{equation}
Note that multiples of $u_{\pi(j)}$, and hence multiples of $\Delta(\lambda)$ for $\lambda\in\pi^{-1}\{\pi(j)\}$, do not change the represented vector. As in the proof of Proposition \ref{cProp:Image of corresponding tropical curve is tropicalization of image of curve} let $r_1< \dots< r_k$ be the elements of
$\{\nu(a_j-a_\lambda)\mid\lambda\in J\}$, and $D_i=\{\lambda\in J\mid \nu(a_j-a_\lambda)\geq r_i\}$ for $1\leq i\leq k$. Furthermore, let $v_i$ be the inclusion-minimal vertex containing $D_i$. The new thing here is that now $r_k=\infty$ is infinite. But the techniques of Proposition \ref{cProp:Image of corresponding tropical curve is tropicalization of image of curve} still work, showing that $\nu(v_i)=r_i$ for $1\leq i<k$, that $v_i\cap J= D_i$ for all $i$, and that there is no vertex $v$ strictly between $v_{i}$ and $v_{i+1}$ such that $v\cap J\neq D_{i+1}$.
With this notation, $\trop(\ev j(\mathcal C))$ is represented by
\begin{equation*}
\sum_{\rho\in\Sigma(1)}\nu(c_\rho)u_\rho +\sum_{i=2}^{k-1}\left(\nu(v_i)-\nu(v_{i-1})\right)s_{v_i},
\end{equation*}
which is easily seen to be equal to $q_{v_{k-1}}$. Because $v_k=\{j\}$, the properties of the $v_i$ imply that $v\cap J= \{j\}$ for every vertex strictly between $v_{k-1}$ and $\{j\}$. This implies that $q_{v_{k-1}}$ differs from the images of the points on $e_j$ only by a multiple of $u_{\pi(j)}$, which finishes the proof.
\end{proof}

Now that we have seen that we have defined the tropical evaluation map $\tev j$ at a label $j\in J$ in the correct way, at least from an algebraic point of view, let us point out its geometric meaning. To do this let us briefly talk about tropical toric varieties. Given a fan $\Theta$ in $\ttorus$, the construction of the toric variety $X_\Theta$ can be mimicked in tropical geometry \cite{P09b}: for every cone $\theta\in\Theta$ one can define an ``affine'' tropical variety $\mathfrak U_\theta=\Hom_{\operatorname{sg}}(\theta^\vee\cap \nt^\vee,\overline{\mathds R})$ and give it the topology induced by the product topology. These topological spaces can be glued, as in the algebraic situation, to a space $\mathfrak X_\Theta$. Again as in toric geometry, there is a natural stratification $\mathfrak X_\Theta=\coprod_\theta \tOrb(\theta)$, where $\tOrb(\theta)=\nt/\Lin(\theta)$ is defined as before. The algebraic and tropical constructions can be connected via the so-called extended tropicalization $\trop:X_\Theta\rightarrow \mathfrak X_\Theta$ which is equal to the tropicalization $O(\theta)\rightarrow \tOrb(\theta)$ as defined above on each stratum. Returning to our situation, the evaluation map $\tev j$ assigns to a curve $(\mathfrak C, (e_\lambda),\tpara)$ the topological limit of the ray $e_j$ in $\mathfrak X_\Sigma$. This limit point lies in the boundary of $\mathfrak X_\Sigma$, in the stratum $\tOrb(j)$ to be precise, and could be considered as the image point of the foot of $e_j$. Thus, the fibers of $\tev j$  consist precisely of the curves whose foot with label $j$ is mapped to a fixed point in the boundary, that is they have the same interpretation as their algebraic counterparts.

The following corollary summarizes our results on rational curves by giving a quite general correspondence theorem. To shorten notation, we define $O(i)\coloneqq T$ and $\tOrb(i)\coloneqq\ttorus$ for $i\in I$.

\begin{corollary}
For every $\lambda\in L_0$ let $A_\lambda$ be a rational affine linear subspace of $\tOrb(\lambda)$ of codimension $d_\lambda$  with difference space $L_\lambda$. Furthermore, assume that 
\begin{equation*}
\sum_{\lambda\in L_0} d_\lambda = |I|+|J|-3 +\dim(\ttorus). 
\end{equation*}
Then for a generic family of points $(p_\lambda)\in \prod _\lambda O(\lambda)$ with $\trop(p_\lambda)\in A_\lambda$ for all $\lambda\in L_0$ we have
\begin{equation*}
\left|\bigcap_{\lambda\in L_0} \ev \lambda^{-1} \left(\mathds G(L_\lambda\cap N_{\tOrb(\lambda)})p_\lambda\right)\right|=\prod_{\lambda\in L_0} \tev \lambda ^*(A_\lambda) \cdot \tMlno I{\ttorus}{\Delta},
\end{equation*}
where $\mathds G(L_\lambda\cap N_{\tOrb(\lambda)})$ denotes the subtorus of $O(\lambda)$ generated by the 1-psgs in $L_\lambda\cap N_{\tOrb(\lambda)}$.
\end{corollary}

\begin{proof}
The proof is analogous to that of Corollary \ref{tCor:Tropical and Algebraic Numbers are Equal}.
\end{proof}

\subsection{Lines on Cubic Surfaces}
It is a well-known theorem by Cayley and Salmon that in algebraic geometry every smooth cubic surface in projective three-space over an algebraically closed field contains exactly $27$ straight lines. The strong analogy between algebraic and tropical geometry makes it promising to prove a similar result in tropical geometry. However, it has already been discovered that the exact analogue of the algebraic statement is simply wrong in the tropical world. It has been shown in \cite{Vig10} that there are smooth tropical cubic surfaces that contain infinitely many lines. On the other hand, in \cite{Vig07} the author shows that general smooth  tropical cubic surfaces of a certain combinatorial type contain exactly $27$ lines. This encourages to believe that it is still possible to prove a slightly modified version of the statement, for example by allowing the lines on tropical cubic surfaces to have multiplicities. These multiplicities should, of course, be $0$ for all but finitely many lines on every cubic surface. Also, they should be defined purely tropically, for example by intersection-theoretical considerations in some moduli space. Finally, it would be nice if the multiplicities would reflect the relative realizability of the lines, that is whenever a tropical cubic $\mathfrak X$ is the tropicalization of a generic algebraic cubic $X$, the multiplicity of each tropical line in $\mathfrak X$ should be equal to the number of algebraic lines in $X$ tropicalizing to it. One way of defining these multiplicities, namely by constructing moduli spaces for tropical lines in given tropical surfaces, was suggested in \cite[Section 3.3]{DennisDiss}. In an example class of cubics, these multiplicities indeed add up to $27$. Unfortunately, this has not been proven for general cubics, and it is difficult to decide whether the lines in question are relatively realizable or not. For this reason, we chose a different approach, which also allows to use our previous results. Namely, we tropicalize the algebraic incidence correspondence of lines and cubics.

But first, let us review some results about tropical hypersurfaces. In doing so, we will also recall the definitions of tropical lines and smooth tropical cubics. Assume we are given a tropical polynomial  $\mathfrak f=\sum_{m\in\mathds Z^n} \mathfrak a_m x^m$, where the coefficients are in $\overline{\mathds R}=\mathds R\cup\{\infty\}$ and only finitely many of them are not equal to $\infty$. The vanishing set $|\mathfrak V_a(\mathfrak f)|$ of a $\mathfrak f$ is equal to all points $\mathfrak x\in\mathds R^n$ such that the minimum
\begin{equation*}
\mathfrak f(\mathfrak x)=\min_{m\in\mathds Z^n} \{\mathfrak a_m+\langle m, \mathfrak x\rangle\}
\end{equation*}
is obtained at least twice. This set has a polyhedral structure whose cells are in dimension-reversing correspondence to the cells of the Newton subdivision of the Newton polytope $\Newt(\mathfrak f)=\conv\{m\mid a_m\neq\infty\}$ \cite[Prop. 3.1.6]{TropBook}. It can be made a tropical variety $\mathfrak V_a(\mathfrak f)$ by assigning to a facet of $|\mathfrak V_a(\mathfrak f)|$ the lattice length of its corresponding edge in the Newton subdivision. This construction is compatible with tropicalization, that is if $f=\sum_{m\in\mathds Z^n} a_m x^m$ is a polynomial in the coordinate ring of $\mathds G_m^n$, and we tropicalize its coefficients to obtain $\mathfrak f =\sum_{m\in\mathds Z^n} \nu(a_m) x^m$,  then we have $\Trop\left(V(f)\right)=\mathfrak V_a(\mathfrak f)$ \cite[Lemma 3.4.6]{TropBook}.
Now given a point $\mathfrak p\in\tT n$ with representative $\mathfrak x\in \mathds R^{n+1}$, and a homogeneous polynomial $\mathfrak f\in \overline{\mathds R}[x_0^{\pm1},\dots,x_n^{\pm1}]$, the value $\mathfrak f(\mathfrak x)$ is not independent of the choice of $\mathfrak x$. Nevertheless, it does not depend on $\mathfrak x$ whether or not the minimum is obtained twice. Thus we get a well-defined set $|\mathfrak V_p(\mathfrak f)|\subset \tT n$, which can be made into a tropical variety $\mathfrak V_p(\mathfrak f)$ in the same way as in the affine case.

Unfortunately, unlike in the algebraic case, there is no bijective correspondence between tropical hypersurfaces in $\tT n$ and homogeneous polynomials in $\overline{\mathds R}[x_0^{\pm1},\dots,x_n^{\pm1}]$ modulo multiplication by a monomial. But for a given hypersurface $\mathfrak X\subseteq \tT n$, the Newton polytopes and subdivisions of the tropical polynomials which have $\mathfrak X$ as their vanishing set are translates of each other. Thus it makes sense to talk of the Newton polytope and subdivision of $\mathfrak X$ and to call $X$ smooth if the Newton subdivision is unimodular.

Following \cite{Vig10} we define tropical cubic surfaces in $\tT 3$ as the vanishing sets of tropical polynomials in the tropical polynomial ring $\overline{\mathds R}[x_0,\dots, x_3]$ whose Newton polytopes are equal to three times the standard simplex, that is equal to $3\conv\{t_0,\dots,t_3\}$, where $(t_i)_i$ denotes the standard basis. 

The second important class of objects relevant to counting lines in cubics is, of course, that of lines in $\tT 3$. Such a line, or more generally a tropical line in $\tT n$ is defined as a genus $0$ curve with exactly $n+1$ unbounded edges that have weight $1$, and as directions the images $s_0,\dots,s_n$ in $\tT n$ of the standard basis vectors of $\mathds R^{n+1}$. Equivalently, a tropical line is the image of a parametrized curve in $\tMlno 0n1$.

Now that we have established which tropical objects we actually want to count, we can start relating them to their analogous algebraic objects. For this we will work over the field $K=\mathds C((t^{\mathds R}))$ of generalized power series. As indicated earlier, we want to tropicalize the algebraic incidence correspondence $M$ of lines and cubic surfaces. This variety parametrizes all pairs $(L,X)$ of lines $L$ and cubic surfaces $X$ such that $L\subseteq X$. As lines in $\PS 3$ are parametrized by the Grassmannian $\G(2,4)$, and cubics are uniquely determined by the $20$ coefficients of the homogeneous polynomials of degree $3$ in $K[x_0,\dots,x_3]$ generating their vanishing ideals, two such polynomials cutting out the same cubic if and only if they differ by an element in $K^*$, the incidence correspondence $M$ is a subvariety of $\G(2,4)\times \PS {19}$. Since only lines whose Plücker coordinates are all nonzero tropicalize to tropical lines in our sense, and only cubics whose defining polynomials contain all monomials of degree $3$ can possibly tropicalize to smooth tropical cubics, we can replace $M$ by the part of the incidence correspondence contained in $\G^0(2,4)\times\aT{19}$. 

Consider $\G^0(2,4)$ in its Plücker embedding in $\aT 5$. By \cite[Thm. 3.8]{SS04}, the tropicalization $\tG(2,4)$ of the Grassmannian $\G^0(2,4)$ parametrizes the tropical lines in a way respecting the tropicalization. The situation is not quite as nice for cubics: A given point $f\in\PS{19}$, representing a cubic polynomial modulo constant factor, tropicalizes to a point $\mathfrak f=\trop(f)\in\tT {19}$, determining a tropical cubic polynomial modulo constant summand. Of course, the vanishing set of a tropical polynomial does not change when adding a constant summand. Thus, as in the algebraic setting, there is a well-defined tropical hypersurface $\mathfrak V_p(\mathfrak f)$, which is clearly equal to $\Trop\left(V_p(f)\right)$. However, as mentioned above, $\tT {19}$ does not parametrize the set of tropical cubics since the map assigning to $\mathfrak f\in \tT{19}$ its associated cubic $\mathfrak V_p(\mathfrak f)$ is not injective. Fortunately, if $\mathfrak X$ is a smooth cubic, there actually is a unique $\mathfrak f\in \tT{19}$ such that $\mathfrak V_p(\mathfrak f)=\mathfrak X$. 

Let us finally apply our results from Section \ref{tSec:Relating Algebraic and Tropical Enumerative Geometry} to the $19$-dimensional subvariety $M$ of $T= \aT 5\times \aT {19}$. We only consider one evaluation map $\ev{}$ (which is not really an evaluation here), namely the restriction to $M$ of the projection $\extev {}=\pr:T\rightarrow \aT {19}$. Of course, $\extev{}$ is dominant and $\aT {19}$ has the correct dimension to apply Corollary \ref{tCor:Tropical and Algebraic Numbers are Equal}. By Cayley and Salmon's theorem, the fiber $\ev{}^{-1}\{f\}$ consists of exactly $27$ points for all $f$ in the open subset of $\aT {19}$ of smooth cubic surfaces. So if we let $\mathfrak M=\Trop(M)$, and let $\tev{}$ be the restriction of $\textev{}=\Trop(\extev{})$ to $\mathfrak M$, Corollary \ref{tCor:Tropical and Algebraic Numbers are Equal} tells us that $\deg(\tev{}^*[\mathfrak f]\cdot\mathfrak M)=27$ for every $\mathfrak f\in\tT{19}$. Every point  $(\mathfrak L,\mathfrak f)\in\mathfrak M\subseteq \tG(2,4)\times\tT{19}$ is the tropicalization of a pair $(L,f)\in M$. In particular, we have $\mathfrak L\subseteq \mathfrak V_p(\mathfrak f)$. We see that for every tropical cubic polynomial $\mathfrak f$ we get multiplicities on the tropical lines contained in $\mathfrak V_p(\mathfrak f)$ that add up to $27$ and respect tropicalization. Since a smooth tropical cubic corresponds to a unique cubic tropical polynomial, we obtain the following proposition:

\begin{proposition}
Let $\mathfrak X$ be a smooth tropical cubic. Then there are intersection-theoretically defined multiplicities $\mult(\mathfrak L)$ for the lines $\mathfrak L\subseteq \mathfrak X$ that are compatible with tropicalization. That is, if $X$ is a generic smooth cubic surface tropicalizing to $\mathfrak X$, then the multiplicity $\mult(\mathfrak L)$ of a tropical line $\mathfrak L\subseteq \mathfrak X$ is equal to the number of lines of $X$ tropicalizing to $\mathfrak L$. In particular, only finitely many tropical lines in $\mathfrak X$ have nonzero weight, the weights are all nonnegative, and they add up to $27$.
\end{proposition}

\vfill
\section*{Index of Notation}
\renewcommand\arraystretch{1.4}
{\small
\begin{tabular}{@{}>{$}l<{$}p{.62\textwidth}p{\widthof{Construction 3.1}}@{}}
\dPair S & set of pairs of distinct elements of $S$ & page \pageref{pSec:Preliminaries} \\
(K,\nu) & ground field with valuation (algebraically closed of characteristic $0$) & page \pageref{pSec:Preliminaries}\\
\VG & value group & page \pageref{pSec:Preliminaries}\\
(R,\mathfrak m),~ \kappa & valuation ring of $K$, residue field & page \pageref{pSec:Preliminaries}\\
\aT r & the algebraic torus $\mathds G^{r+1}_m/\mathds G_m$ & Subsection \ref{subSec:Tropical Varieties}\\
\tT r & the tropical torus $\mathds R^{r+1}/\mathds R \mathbf 1$ & Subsection \ref{subSec:Tropical Varieties}\\
\VG(\ttorus) & $\VG$-rational points in $\ttorus$ & Subsection \ref{pSubsec:Tropicalization}\\
\ini_{\mathfrak t}X & initial degeneration of $X$ at $\mathfrak t$ & Subsection \ref{pSubsec:Tropicalization}\\
\mathcal T & standard $R$-model of an algebraic torus $T$ & Subsection \ref{pSubsec:Tropicalization}\\
\tMlI {I}  & moduli space of labeled parametrized $I$-marked rational tropical curves in $\tT r$ of degree $d$ & Subsection \ref{cSubSec:Tropical Curves}\\
\Mlno I & moduli space of labeled parametrized $I$-marked rational algebraic curves in $\PS r$ of degree $d$ & Subsection \ref{cSubSec:Algebraic Cuves} \\
\tMlno I{\ttorus}{\Delta} &  moduli space of labeled parametrized $I$-marked rational tropical curves in the tropical torus $\ttorus$ of degree $\Delta$ & Subsection \ref{cSubSec:Tropical Curves}\\
\MlIrd I{X}c & moduli space of labeled parametrized $I$-marked rational algebraic curves in the toric variety $X$ satisfying the tangency condition $c$ & Subsection \ref{aSubSec:GeneralDegrees}\\
\ev i ,\tev i & algebraic/tropical evaluation at the $i$-th marked point  \newline $i$-th abstract evalution map & Sections \ref{cSec:Curves}, \ref{aSec:Applications} \newline  Section \ref{tSec:Relating Algebraic and Tropical Enumerative Geometry}\\
 \extev i  & extended algebraic $i$-th evaluation map \newline $i$-th extended abstract evaluation map & pages \pageref{cProp:Evaluation maps can be extended}, \pageref{aProp:Evaluation maps can be extended} \newline Section \ref{tSec:Relating Algebraic and Tropical Enumerative Geometry}\\
 \textev i & extended tropical $i$-th evaluation map $\textev i=\Trop(\extev i)$ & Section \ref{tSec:Relating Algebraic and Tropical Enumerative Geometry}\\
\tpl & tropical Plücker embedding & Subsection \ref{cSubSec:Tropical Curves}\\
\mathfrak C_{\mathcal C} & corresponding tropical curve to the algebraic curve $\mathcal C$ & Construction \ref{cConstr:Corresponding tropical curve}\\
\end{tabular}
\renewcommand\arraystretch{1.0}
}

\bibliography{}

\providecommand{\bysame}{\leavevmode\hbox to3em{\hrulefill}\thinspace}
\providecommand{\MR}{\relax\ifhmode\unskip\space\fi MR }
\providecommand{\MRhref}[2]{%
  \href{http://www.ams.org/mathscinet-getitem?mr=#1}{#2}
}
\providecommand{\href}[2]{#2}
\begin{thebibliography}{GKM09}

\bibitem[AHR16]{AHR16}
L.~{Allermann}, S.~{Hampe}, and J.~{Rau}, \emph{{On rational equivalence in
  tropical geometry.}}, {Can. J. Math.} \textbf{68} (2016), no.~2, 241--257.

\bibitem[AR10]{AR10}
L.~Allermann and J.~Rau, \emph{First steps in tropical intersection theory},
  Mathematische Zeitschrift \textbf{264} (2010), no.~3, 633--670.

\bibitem[BG84]{BG84}
R.~Bieri and J.~R.~J. Groves, \emph{The geometry of the set of characters
  induced by valuations}, J. Reine Angew. Math. \textbf{347} (1984), 168--195.

\bibitem[CLS11]{CLS}
D.~Cox, J.~Little, and H.~Schenck, \emph{Toric varieties}, Amer. Math. Soc.,
  Providence, RI, 2011.

\bibitem[Cox95]{Cox95}
D.~A. Cox, \emph{The functor of a smooth toric variety}, Tohoku Math. J. (2)
  \textbf{47} (1995), no.~2, 251--262.

\bibitem[Dra08]{Drai08}
J.~Draisma, \emph{A tropical approach to secant dimensions}, Journal of Pure
  and Applied Algebra \textbf{212} (2008), no.~2, 349--363.

\bibitem[FR13]{FR10}
G.~{Fran\c{c}ois} and J.~{Rau}, \emph{{The diagonal of tropical matroid
  varieties and cycle intersections.}}, {Collect. Math.} \textbf{64} (2013),
  no.~2, 185--210.

\bibitem[Fra13]{F11}
G.~Fran\c{c}ois, \emph{{Cocycles on tropical varieties via piecewise
  polynomials.}}, Proc. Am. Math. Soc. \textbf{141} (2013), no.~2, 481--497.

\bibitem[Ful98]{F98}
W.~Fulton, \emph{Intersection {T}heory}, second ed., Springer, Berlin, 1998.

\bibitem[GKM09]{GKM09}
A.~Gathmann, M.~Kerber, and H.~Markwig, \emph{Tropical fans and the moduli
  spaces of tropical curves}, Compos. Math. \textbf{145} (2009), no.~1,
  173--195.

\bibitem[GM82]{GM82}
I.~M. Gelfand and R.~D. MacPherson, \emph{Geometry in {G}rassmannians and a
  generalization of the dilogarithm}, Adv. in Math. \textbf{44} (1982), no.~3,
  279--312.

\bibitem[GM08]{GM08}
A.~Gathmann and H.~Markwig, \emph{Kontsevich's formula and the {WDVV} equations
  in tropical geometry}, Adv. Math. \textbf{217} (2008), no.~2, 537--560.

\bibitem[Gub13]{Gub12}
W.~Gubler, \emph{A guide to tropicalizations}, Algebraic and Combinatorial
  Aspects of Tropical Geometry, Contemporary Mathematics, vol. 589, Amer. Math.
  Soc., Providence, RI, 2013, pp.~125--189.

\bibitem[Mik05]{Mik05}
G.~Mikhalkin, \emph{Enumerative tropical algebraic geometry in
  {$\mathbb{R}^2$}}, J. Amer. Math. Soc. \textbf{18} (2005), no.~2, 313--377.

\bibitem[MS15]{TropBook}
D.~Maclagan and B.~Sturmfels, \emph{Introduction to {T}ropical {G}eometry},
  Graduate Studies in Mathematics, vol. 161, American Mathematical Society,
  Providence, RI, 2015.

\bibitem[NS06]{NS06}
T.~Nishinou and B.~Siebert, \emph{Toric degenerations of toric varieties and
  tropical curves}, Duke Math. J. \textbf{135} (2006), no.~1, 1--51.

\bibitem[Och13]{DennisDiss}
D.~Ochse, \emph{Moduli spaces of rational tropical stable maps into smooth
  tropical varieties}, Ph.D. thesis, TU Kaiserslautern, 2013,
  \url{https://kluedo.ub.uni-kl.de/files/3528/thesis+dennis+ochse-1.pdf}.

\bibitem[OP13]{OP13}
B.~Osserman and S.~Payne, \emph{Lifting tropical intersections}, Doc. Math.
  \textbf{18} (2013), 121--175.

\bibitem[OR13]{OR11}
B.~Osserman and J.~Rabinoff, \emph{Lifting non-proper tropical intersections},
  {Tropical and Non-Archimedean Geometry}, Contemporary Mathematics, vol. 605,
  Amer. Math. Soc., Providence, RI, 2013, pp.~15--44.

\bibitem[Pay09a]{P09b}
S.~Payne, \emph{Analytification is the limit of all tropicalizations}, Math.
  Res. Lett. \textbf{16} (2009), no.~3, 543--556.

\bibitem[Pay09b]{P09}
\bysame, \emph{Fibers of tropicalization}, Math. Z. \textbf{262} (2009), no.~2,
  301--311.

\bibitem[Sha13]{Sh13}
K.~Shaw, \emph{{A tropical intersection product in matroidal fans.}}, {SIAM J.
  Discrete Math.} \textbf{27} (2013), no.~1, 459--491.

\bibitem[Shu05]{Shustin05}
E.~Shustin, \emph{A tropical approach to enumerative geometry}, Algebra i
  Analiz \textbf{17} (2005), no.~2, 170--214.

\bibitem[Spe05]{Speyer-diss}
D.~E. Speyer, \emph{Tropical geometry}, Ph.D. thesis, UC Berkeley, 2005.

\bibitem[SS04]{SS04}
D.~E. Speyer and B.~Sturmfels, \emph{The tropical {G}rassmannian}, Adv. Geom.
  \textbf{4} (2004), no.~3, 389--411.

\bibitem[ST08]{ST08}
B.~{Sturmfels} and J.~{Tevelev}, \emph{Elimination theory for tropical
  varieties.}, Math. Res. Lett. \textbf{15} (2008), no.~3, 543--562.

\bibitem[Tyo12]{Tyomkin12}
I.~Tyomkin, \emph{Tropical geometry and correspondence theorems via toric
  stacks}, Math. Ann. \textbf{353} (2012), no.~3, 945--995.

\bibitem[Vig07]{Vig07}
M.~D. Vigeland, \emph{Tropical lines on smooth tropical surfaces}, 2007, \href
  {http://arxiv.org/abs/0708.3847} {\path{arXiv:0708.3847}}.

\bibitem[Vig10]{Vig10}
\bysame, \emph{Smooth tropical surfaces with infinitely many tropical lines},
  Ark. Mat. \textbf{48} (2010), no.~1, 177--206.

\end{thebibliography}

\normalsize
Andreas Gross, Fachbereich Mathematik, Technische Universität Kaiserslautern, Postfach 3049, 67653 Kaiserslautern, Germany \texttt{agross@mathematik.uni-kl.de}

\end{document}